\documentclass[]{elsarticle}

\usepackage{lineno,hyperref}

\usepackage{graphicx}
\usepackage{bm}
\usepackage{multirow}
\usepackage{amsmath}
\usepackage{amsfonts}
\usepackage{amssymb}
\usepackage{amsthm}
\usepackage{secdot}

\journal{Computer-Aided Design}

\theoremstyle{plain}

\theoremstyle{definition}

\theoremstyle{remark}

\numberwithin{equation}{section}
\numberwithin{theorem}{section}
\numberwithin{remark}{section}

\begin{document}

\begin{frontmatter}

\title{Surface and Hypersurface Meshing Techniques for Space-Time Finite Element Methods}

\author{Jude T. Anderson \corref{mycorrespondingauthor}}

\cortext[mycorrespondingauthor]{Corresponding author}
\ead{jta29@psu.edu}

\address{Department of Mechanical Engineering, The Pennsylvania State University, University Park, Pennsylvania 16802}

\author{David M. Williams} 

\address{Department of Mechanical Engineering, The Pennsylvania State University, University Park, Pennsylvania 16802}

\author{Andrew Corrigan}

\address{Laboratories for Computational Physics and Fluid Dynamics, Naval Research Laboratory, Washington, DC 20375}
\fntext[fn1]{Distribution Statement A: Approved for public release. Distribution is unlimited.}

\begin{abstract}
A general method is introduced for constructing two-dimensional (2D) surface meshes embedded in three-dimensional (3D) space time, \emph{and} 3D hypersurface meshes embedded in four-dimensional (4D) space time. In particular, we begin by dividing the space-time domain into time slabs. Each time slab is equipped with an initial plane (hyperplane), in conjunction with an unstructured simplicial surface (hypersurface) mesh that covers the initial plane. We then obtain the vertices of the terminating plane (hyperplane) of the time slab from the vertices of the initial plane using a space-time trajectory-tracking approach. Next, these vertices are used to create an unstructured simplicial mesh on the terminating plane (hyperplane). Thereafter, the initial and terminating boundary vertices are stitched together to form simplicial meshes on the intermediate surfaces or \emph{sides} of the time slab. 
After describing this new mesh-generation method in rigorous detail, we provide the results of multiple numerical experiments which demonstrate its validity and flexibility.
\end{abstract}


\end{frontmatter}


\section{Introduction}
\label{sec;introduction}

Since its inception, the finite element method has often been limited to stationary three-dimensional (3D) geometries due to the available meshing capabilities.  
However, in the last two decades, research has been conducted with the goal of accurately simulating fluid-structure interactions (FSI) for 3D moving bodies. Towards this end,  one may extrude or extend a 3D object along the temporal direction in order to capture its movement in a four-dimensional (4D) space-time setting.  As one may imagine, this process is not intuitive, as the entirety of the domain is no longer directly visible and can only be observed through projections or hyperplane cross sections.  Furthermore, extending the current technology to properly mesh these domains has proven to be a difficult task. The existing meshing technologies include methods for generating structured and semi-unstructured 4D meshes; however, the literature does not appear to contain a method for fully-unstructured mesh generation with boundary recovery in 4D. In what follows, we briefly review some of the relevant work on space-time volume meshing and classical surface meshing; then we provide an overview of our current efforts to extend this work to create fully-unstructured, 4D meshes. 

Some of the earliest work related to space-time finite element methods in one spatial dimension plus time (1D+$t$) can be found in the papers of Hughes and Hulbert~\cite{hughes1988space,hulbert1990space}. Thereafter, Behr~\cite{behr2008simplex} developed a method for semi-unstructured temporal extrusion that applies to both two-dimensional (2D) and 3D meshes. Broadly speaking, Behr introduced a process for extruding the triangular elements of a 2D surface mesh along the temporal direction to create triangular prisms that can each be discretized into tetrahedra conforming to the Delaunay criterion.  The process is similar for 3D hypersurface meshes, where a 3D hypersurface mesh of tetrahedra are extruded along the temporal direction to form 4D tetrahedral prism elements, which are subsequently discretized into pentatopes also conforming to the Delaunay criterion. This approach has been successfully applied to a wide range of applications, as evidenced by the work in~\cite{pauli2017stabilized,von2019simplex,karyofylli2019simplex,make2022spline,von2022time}. 

The prism extrusion method of Behr was expanded upon by von Danwitz et al.~\cite{von2021four}. In particular, they extend the method to accommodate time-variant topology through what they call a 4D-elastic mesh update method.  Essentially, this does not change the connectivity of the 4D mesh but merely deforms the existing elements to conform to the varying surface topology. The most recent work on this particular topic appears to be that of Karyofylli and Behr \cite{karyofylli2022simplex}.

In addition, a number of researchers have extended the extrusion-based method to accommodate rotational motions. For example, Wang and Persson~\cite{wang2015high} employ a similar method to Behr in 2D+$t$ in order to generate an initial tetrahedral mesh; thereafter, they subdivide the mesh into a stationary region, a rotating region, and a buffer region that resides at the interface between the two. During rotation, the connectivity between the rotational region and the stationary region is maintained via reconnections (edge flips or face flips) in the buffer region. Wang and Persson's approach is essentially a space-time, \emph{sliding-mesh} approach. It has been extended to 3D+$t$ for very simple cases~\cite{wang2015discontinuousthesis}. Its viability appears to hold only for applications in which the boundary motion is purely rotational, and is known \emph{a priori}. We note that a very similar sliding-mesh approach has been recently developed by Horváth and Rhebergen~\cite{horvath2022conforming}. This work appears to extend, and in some ways improve upon the previous work of Wang and Persson.

In contrast to the extrusion-based methods (above), Foteinos and Chrisochoides~\cite{foteinos20154d} were able to generate unstructured 4D hypervolume meshes using a typical Delaunay-based mesh generator up-scaled to accommodate four dimensions.  Although this work is significant, it does not present a clear mechanism for recovering the boundary, i.e.~recovering the surface mesh that resides on the boundary. This issue of `boundary recovery' is a common problem in Delaunay-based meshing techniques.  Many researchers, such as Si et al.~\cite{si20113d} and Liu et al.~\cite{liu2014boundary}, detail various strategies for recovering a surface mesh in 3D, in conjunction with a Delaunay mesh generator. However, due to the lack of boundary recovery strategies in 4D, extrusion-based methods similar to Behr's remain dominant.

As another alternative to the methods proposed above, traditional advancing front techniques~\cite{lohner1988generation, george1994advancing, lo2014finite} were expanded upon to create a space-time mesh generation method that is known as \emph{pitching tents}~\cite{ungor2002pitching, erickson2005building, abedi2004spacetime}. Here, each vertex in the spatial domain is projected along the temporal direction to generate a new vertex.  New elements (one dimension higher than the original mesh) are created by tessellating the region formed by the neighboring faces of the original vertex and the new vertex. Recently, this method has been applied to hyperbolic systems~\cite{gopalakrishnan2015tent,gopalakrishnan2017mapped,drake2022convergence} and the Maxwell equations~\cite{gopalakrishnan2020explicit}. These methods are usually best-suited for wave-propagation problems.


Most of the existing research (above) focuses on the generation of space-time volume meshes in 2D+$t$ and hypervolume meshes in 3D+$t$. To our knowledge, researchers have not rigorously explored techniques for generating space-time surface meshes in 2D+$t$ and hypersurface meshes in 3D+$t$. Part of the reason for this omission is that boundary conformity is automatically enforced for extrusion-based meshing approaches for problems where the boundary is stationary. Furthermore, boundary motion can (sometimes) be accommodated via the aforementioned elasticity-based approach.  Of course, this comes at the cost of failing to accommodate arbitrary, large-scale boundary motions. More importantly, most volume or hypervolume meshes generated via extrusion are not fully unstructured in both space and time. Therefore, there is some incentive to investigate general surface meshing techniques which can accommodate fully-unstructured, constrained-Delaunay meshing strategies in 4D. 

Naturally, there is already a wealth of literature which discusses surface meshing techniques for traditional, stationary, 3D applications. The majority of the work in this area relies on a method known as parametric mapping, which involves constructing a surface mesh for a 3D application on a reference 2D domain according to a specified metric.  Once the 2D domain is meshed, it is mapped to the 3D domain.  Interesting applications and discussions of this method can be found in~\cite{borouchaki2000parametric, borouchaki1996unstructured, tristano1998advancing, zheng1996three, lee1998automatic, canann1997automatic, cuilliere1998adaptive}. Variations and improvements to this method include separating a surface into patches~\cite{lee2003automatic}, using high-order elements~\cite{sherwin2002mesh}, and using Voronoi diagrams~\cite{levy2013variational}, among other notable works~\cite{borouchaki2005surface, frey2003surface, zhong2014anisotropic}. In addition, Lan and Lo~\cite{lan1996finite} and Cass et al.~\cite{cass1996generalized} developed their own alternatives to parametric mapping that remain in 3D space and employ techniques such as curvature sizing functions to generate valid surface meshes.

The key issue with this existing surface meshing literature is that it is generally limited to 2D surface meshes which are embedded in 3D space. Furthermore, the meshing techniques often depend on a detailed knowledge of the underlying CAD, which is easily available for 3D problems, but may not be fully characterized for 4D space-time problems. 

In this work, we discuss a new approach to 4D meshing, specifically the generation of a hypersurface mesh embedded in 3D+$t$ space time. A key component of this method, is that vertices from the previous time slab change their positions in accordance with tracking space-time trajectories, (computed based on the local hypersurface velocity). In principle, the method can successfully track the movement of any 3D object in question. 
In addition, it allows us to create hyperplanes that are stitched together by tetrahedra in order to form a complete, conforming hypersurface mesh on each time slab. Once this hypersurface mesh is generated, we can create a fully-unstructured hypervolume mesh that conforms to the hypersurface mesh on the slab. Note that this latter step will be reserved for subsequent work.  In what follows, we discuss preliminary concepts relating to surface meshing, then move on to 2D+$t$ and 3D+$t$ illustrations of our technique.  We then present the results of some numerical experiments and conclude by suggesting future work.


\section{Preliminaries}

\begin{figure}[h!]
\centering
\includegraphics[width=10cm]{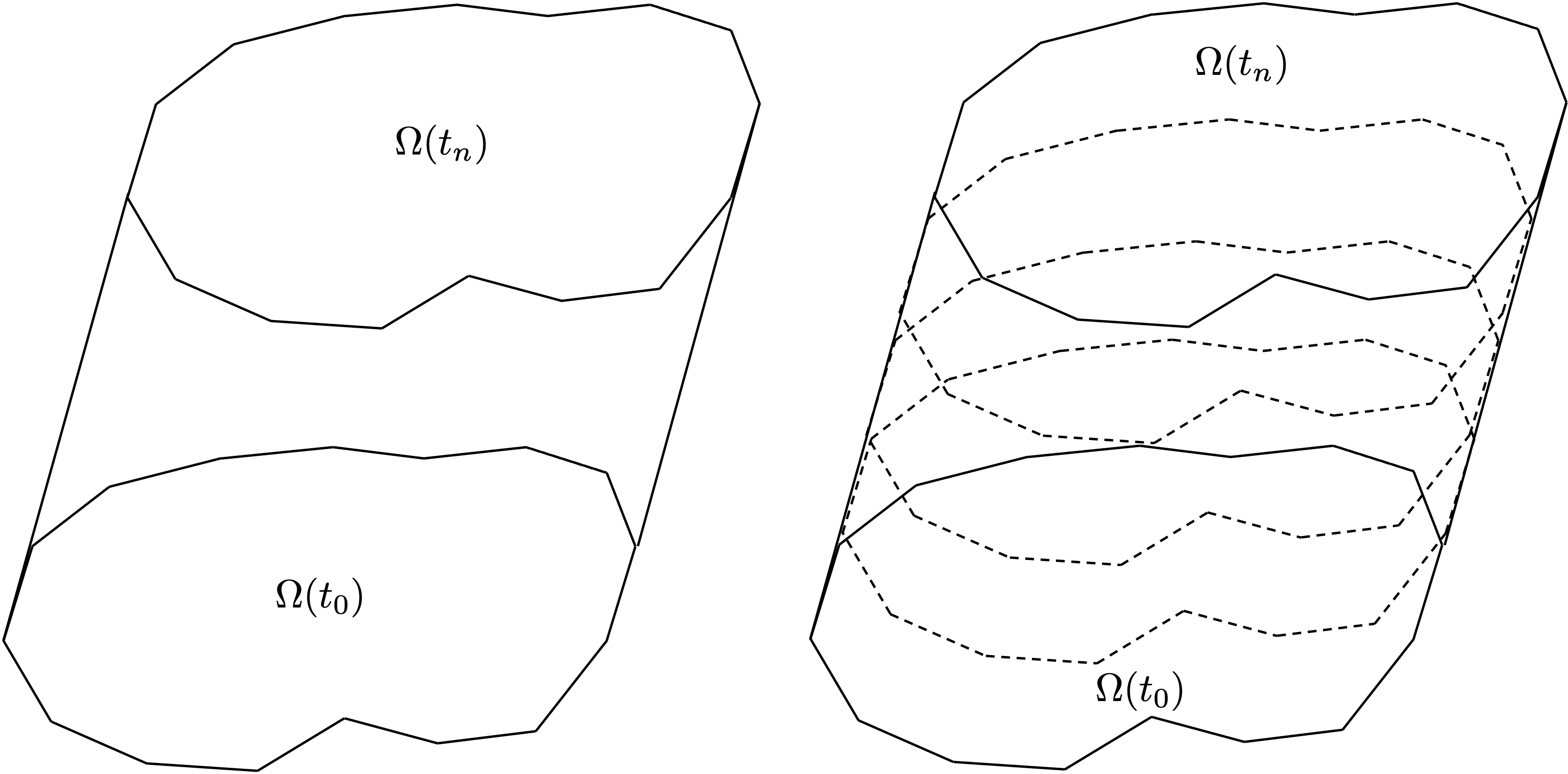}
\caption{Entire space-time domain in 2D+$t$ (left) and subdivision of this domain into space-time slabs (right).}
\label{space_time_slabs}
\end{figure}

We begin by partitioning the space-time domain into slabs. This decomposition process is performed by intersecting the domain with spatial hyperplanes located at regular intervals, (see Figure~\ref{space_time_slabs}). In 3D, each slab contains an ``initial plane" (at $t = t_n$), a ``terminating plane" (at $t = t_{n+1}$), and an ``intermediate surface" which connects the initial plane to the terminating plane. We note that the intermediate surface does not need to be planar. When taken together, the initial plane, terminating plane, and intermediate surface form the boundaries of the space-time slab. These boundaries are 2D surfaces embedded in a 3D space time (2D+$t$). We note that the space-time slabs are generally not formed all at once. Instead, they are formed sequentially on an ``as-needed basis", starting with those at the earliest times, and then continuing on with those at later times. This is a particularly important point, when we consider that the geometry of the space-time domain may not be known \emph{a priori} for certain FSI applications, and a predictor-corrector approach may be necessary to form the topology of the individual space-time slabs, as the simulation evolves.

Before proceeding further, it is important for us to distinguish between ``continuous space-time surfaces" and ``discrete space-time surfaces". A continuous space-time surface is the continuous, CAD definition of a surface, which can only be generated if suitable knowledge of the boundary motion is available. A discrete space-time surface is the discrete, surface mesh that is (frequently) associated with an underlying continuous space-time surface. For the case of 2D+$t$, this surface mesh usually consists of triangles and their associated vertices embedded in 3D space time. For the case of 3D+$t$, the surface is actually a hypersurface which usually consists of tetrahedra and their associated vertices embedded in 4D space time. In this work, we are primarily interested in generating suitable surface meshes (i.e., discrete space-time surfaces). We can summarize our objectives in the following problem statement for 2D+$t$:

\vspace{0.25cm}

\noindent \emph{``Given a surface mesh on the previous space-time slab, find new surface meshes on the initial, intermediate, and terminating surfaces of the next space-time slab, while limiting the amount of space-time CAD information required."}

\vspace{0.25cm}

\noindent We can obtain an equivalent statement for the case of 3D+$t$ by replacing the word ``surface" with the word ``hypersurface" in the statement above. 

\section{Surface and Hypersurface Meshing} \label{surface_mesh_approach}

For the 2D case, we begin by extracting the terminating plane of the previous space-time slab, which is located at $t = t_n$. We assume that this terminating plane is covered with a discrete triangular surface mesh. Of course, if we consider the first space-time slab in our entire space-time domain, the previous space-time slab does not exist. In this case, we simply assume that there is a ghost slab that spans the space from $t = t_{-1}$ to $t = t_0$ and provides us with a terminating surface mesh located at $t = t_0$. Once the terminating surface mesh is identified, our objective is to create new surface meshes on the initial, intermediate, and terminating surfaces of the next space-time slab from $t = t_n$ to $t = t_{n+1}$. With this in mind, we start by setting the surface mesh on the initial plane of our new space-time slab to be identical to the terminating surface mesh of the previous space-time slab. This ensures that the subsequently generated volume mesh on the new space-time slab will maintain conformity with the volume mesh on the previous space-time slab. Next, it is possible to generate the intermediate surface mesh on the ``sides" of the space-time slab, and thereafter, the surface mesh on the terminating plane. In what follows, we will describe the remainder of the surface meshing process in 2D. Thereafter, we discuss the extension to~3D.

\subsection{The Two-Dimensional Case (2D+t)}

In order to begin building the intermediate 2D surface mesh, we extract and mark the edges which represent the discrete boundaries of the surface mesh on the initial plane. Next, we compute the space-time trajectories of these vertices using the local velocity of the space-time CAD surface (which should be known or predicted \emph{a priori}), in conjunction with the well-known ordinary differential equation
\begin{align*}
    \bm{v}\left(t\right) = \frac{d \bm{x}\left(t\right)}{dt}.
\end{align*}
In order to solve this equation, the time interval $dt = t_{n+1}-t_{n}$ is subdivided into $M$ subintervals, where the time-step for each subinterval is simply $dt/M$. On each subinterval, we solve the differential equation above using the latest information about the surface velocity, in conjunction with a standard explicit time-stepping method, such as the forward Euler time-stepping method. In this way, the trajectories of the edge vertices are computed until their location at the final time ($t_{n+1}$) is determined. The final location of the vertices may be impacted by time-integration errors, and therefore, we perform a simple projection procedure to ensure that the vertices lie exactly on the CAD surface at the final time. Note: throughout this process, we assume that the connectivity of the edge vertices does not change. After the vertex trajectories and final locations have been computed, we have two sets of vertices: one on the initial plane at $t = t_n$, and one on the terminating plane at $t = t_{n+1}$. The vertices on the terminating plane are then connected to one another in order to form edges. Thereafter, these edges are connected to the corresponding edges on the initial plane in order to form quadrilateral elements on the intermediate surface. These quadrilateral elements can be subdivided into triangular elements by inserting a Steiner point at the centroid of each quadrilateral element and connecting each Steiner point to the quadrilateral element’s vertices. By following this procedure, we succeed in forming a linearly interpolated surface mesh of triangles on the intermediate surface. Lastly, we collect the edges and vertices on the terminating surface, and send them to a 2D constrained Delaunay mesh generation program (such as Shewchuk’s Triangle program~\cite{shewchuk1996triangle}) in order to generate a surface mesh for the terminating plane. We conclude by synchronizing the connectivity of the initial surface mesh, the intermediate surface mesh, and the terminating surface mesh. The resulting agglomeration of surface meshes provides a hull of triangular elements on the space-time slab from $t = t_n$ to $t = t_{n+1}$.

\begin{figure}[h!]
\centering
\includegraphics[width=13cm]{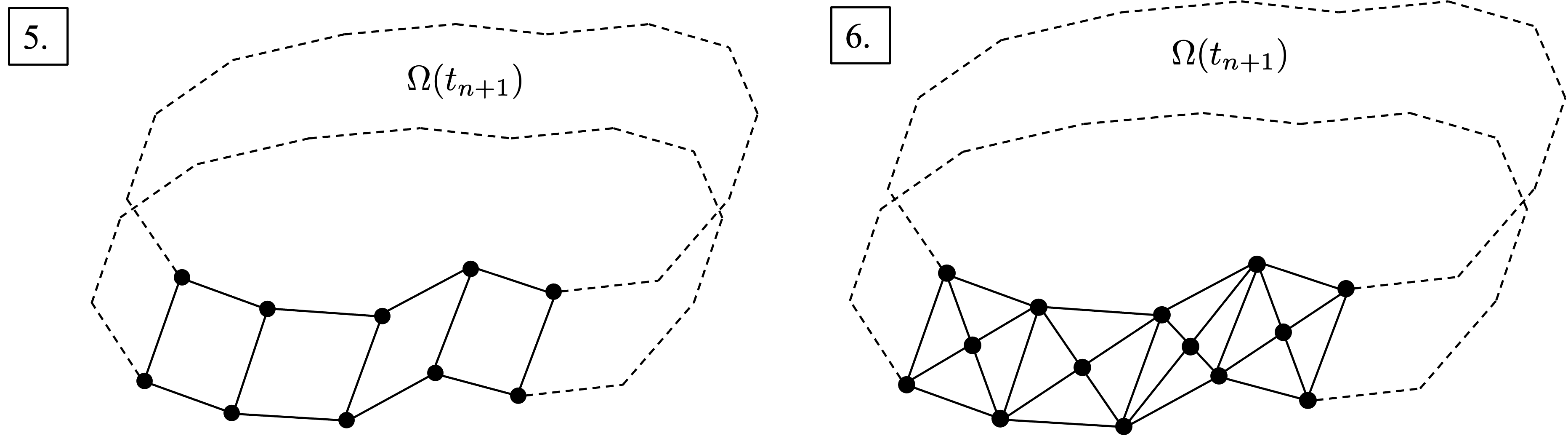}

\vspace{0.5cm}

\includegraphics[width=6.5cm]{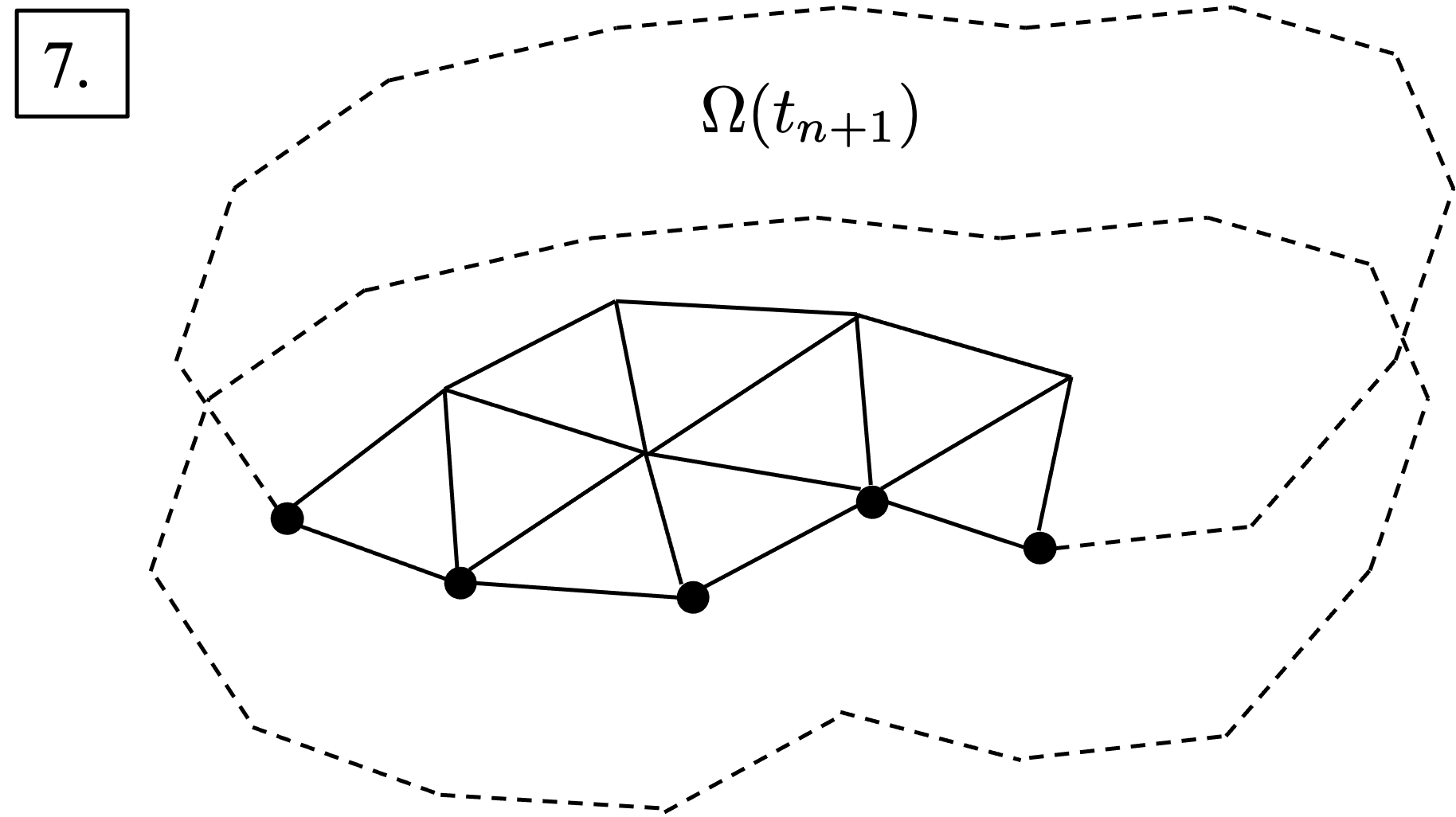}
\caption{An illustration of the process for generating a surface mesh on a space-time slab in 2D+$t$. The numbered steps are explained in the text.}
\label{2D_surface_process}
\end{figure}

The process for generating the surface mesh on a space-time slab is summarized below:

\begin{enumerate}
\item Extract the surface mesh from the terminating plane of the previous space-time slab.
\item Construct the surface mesh for the initial plane of the new space-time slab using the surface mesh from step 1.
\item Extract the boundary edges and vertices of the surface mesh from step 2. Compute the vertex trajectories from $t = t_n$ to $t = t_{n+1}$. Project final point locations to the CAD surface.
\item Connect vertices on the terminating plane to create edges.
\item Connect edges on the terminating plane to edges on the initial plane to create quadrilaterals.
\item Subdivide the quadrilaterals into triangles to generate a triangular surface mesh on the intermediate surface.
\item Use the edges on the terminating plane to generate a triangular surface mesh on the terminating plane.
\end{enumerate}
Thereafter, the surface meshes from steps 2, 6, and 7 are combined to generate a surface mesh for the entire space-time slab. The overall process is shown in Figure~\ref{2D_surface_process}.

\begin{figure}[h!]
\centering
\includegraphics[width=13cm]{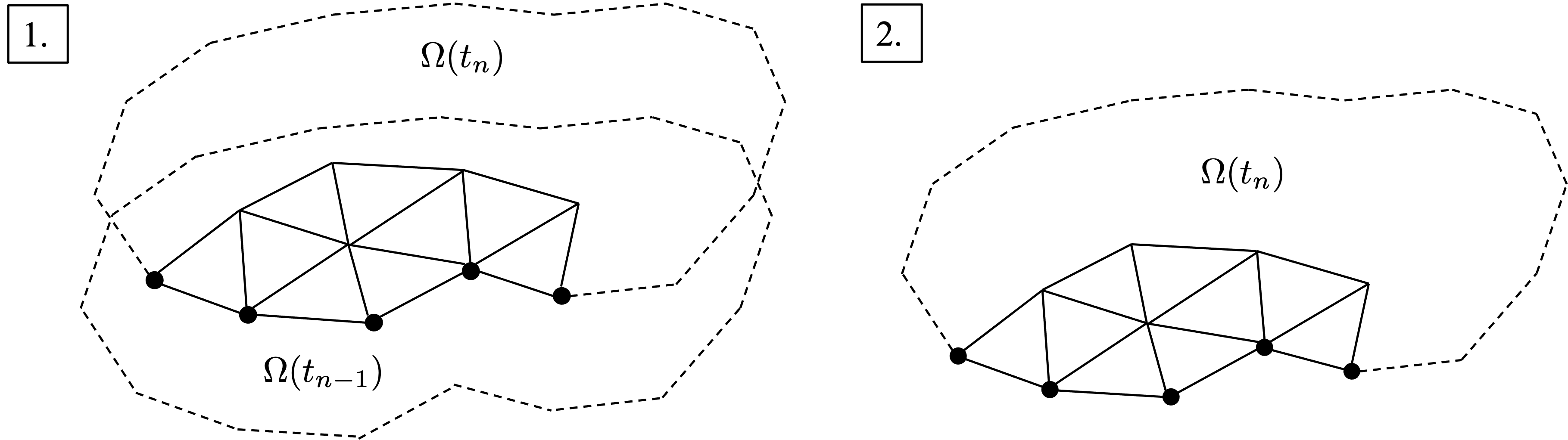}

\vspace{0.5cm}

\includegraphics[width=13cm]{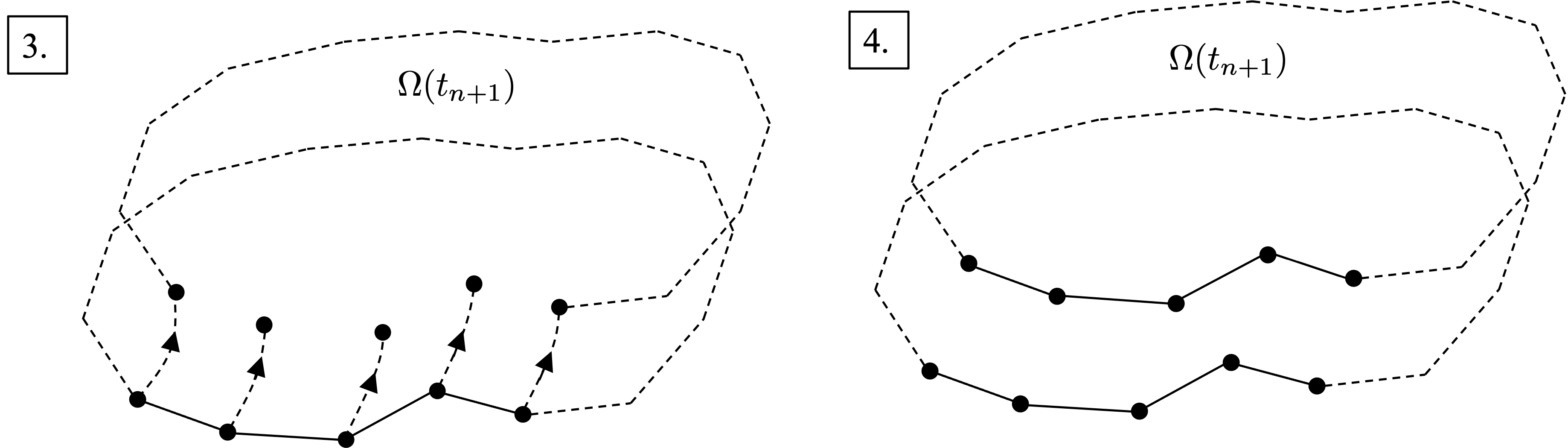}
\end{figure}

\pagebreak

\subsection{The Three-Dimensional Case (3D+$t$)}

\begin{figure}[h!]
\centering
\includegraphics[width=8cm]{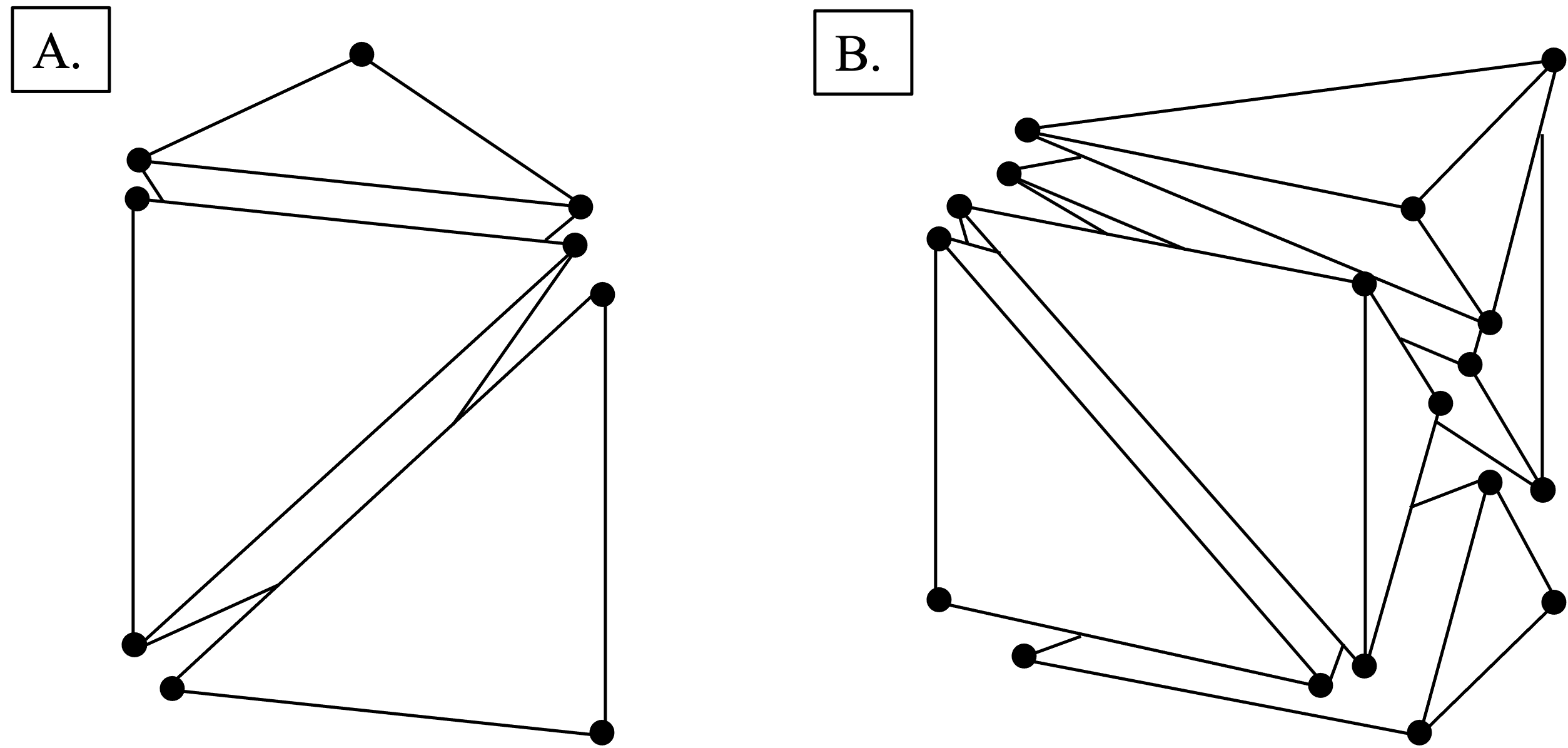}

\vspace{1cm}

\includegraphics[width=8cm]{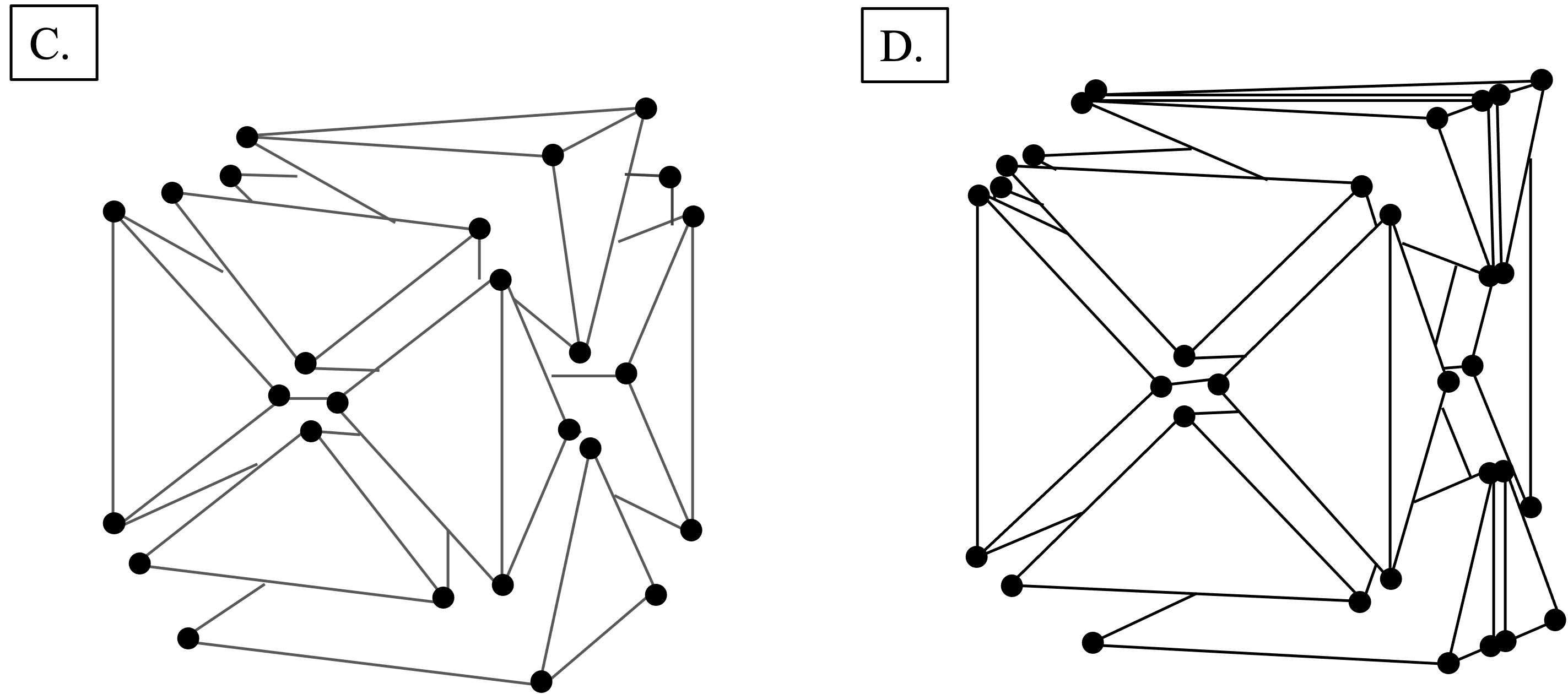}

\vspace{1cm}

\includegraphics[width=3.5cm]{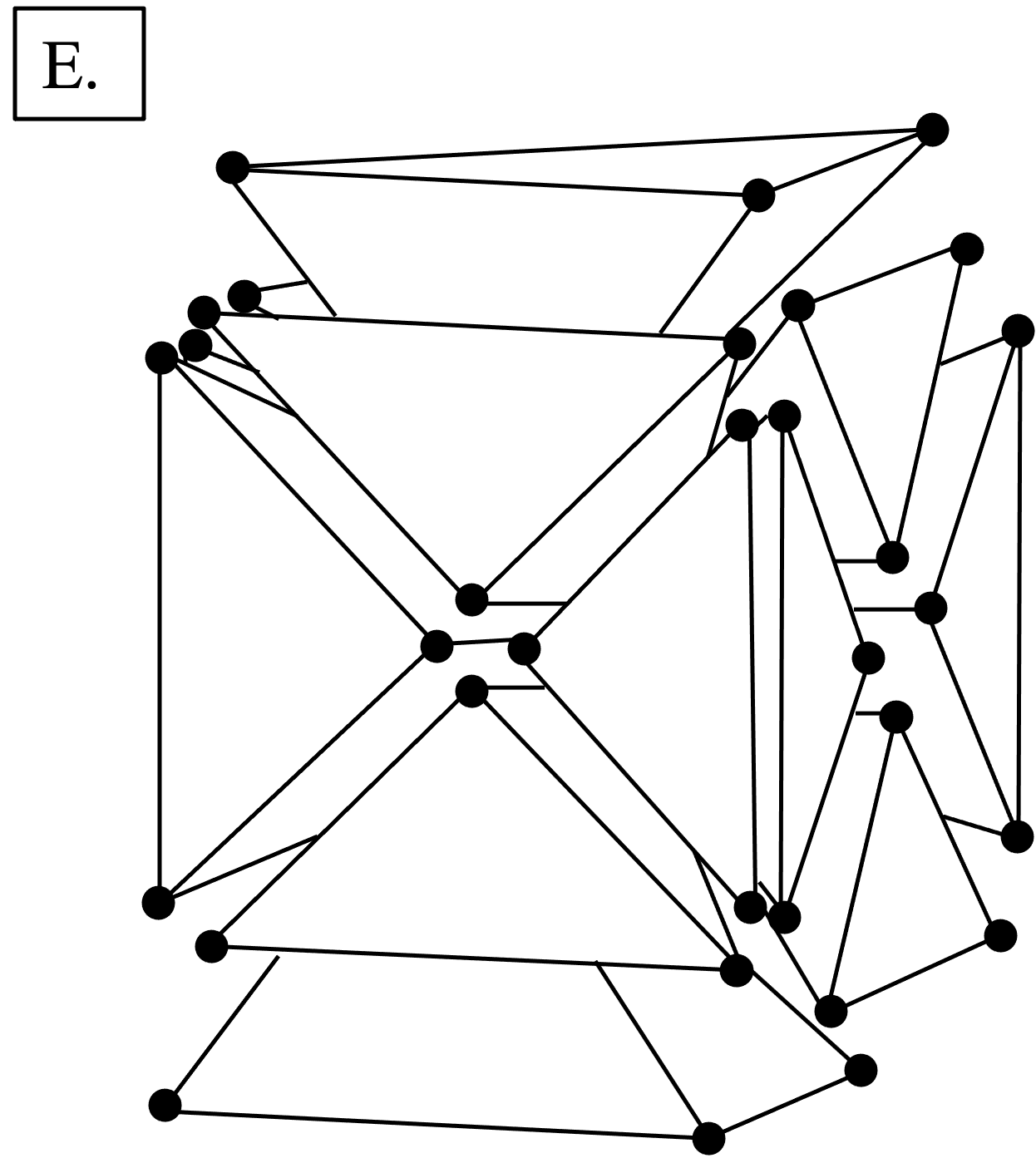}
\caption{Illustrations of the five different strategies for subdividing a triangular prism into tetrahedra elements. The strategies (A-E) subdivide the prism into 3, 6, 10, 12, and 14 tetrahedral elements, respectively.}
\label{prism_splitting}
\end{figure}

In order to build the 3D hypersurface mesh, we first extract the tetrahedral hypersurface mesh on the terminating hyperplane of the previous space-time slab. Following the approach used in the 2D case, we use this hypersurface mesh in order to tessellate the initial hyperplane of the next space-time slab. Thereafter, it remains for us to construct the intermediate hypersurface mesh, and the terminating hypersurface mesh for the space-time slab. Towards this end, we identify the outer boundaries of the initial hypersurface mesh. These boundaries correspond to the set of triangles which lie on the boundaries of the spatial domain at $t = t_n$. Once these triangles have been identified, we can extract their vertices, compute the corresponding space-time trajectories, and project the final point locations to the CAD surface (see the 2D procedure for details). Thereafter, we will have triangle vertices on the initial hyperplane (at $t = t_n$) and on the terminating hyperplane (at $t = t_{n+1}$). The vertices on the terminating hyperplane can be connected to form a triangulation, then the triangles on the initial and terminating hyperplanes can be connected in order to form triangular prisms. Note that these are 3D triangular prisms which are embedded in 4D space time. Once the prisms have been formed, they can be split into tetrahedral elements.  We are aware of at least five different splitting strategies (see Figure~\ref{prism_splitting}). However, an arbitrary splitting is not possible, as it is important to preserve the conformity of adjacent triangular prism faces. Therefore, we require that all splittings of the quadrilateral faces of the triangular prisms are identical under reflections and rotations of the prism unto itself. With this in mind, we prefer two particular splitting techniques. The first technique involves splitting the quadrilateral faces of the triangular prisms into triangles by inserting Steiner points at the centroids of the quadrilateral faces and connecting these points to the quadrilateral’s vertices. Thereafter, the triangular faces of the split prism can be connected to an additional Steiner point located at the prism’s centroid. This results in a total of fourteen tetrahedral elements (see Figure~\ref{prism_splitting}, E). This splitting is natural because it is merely a generalization of the splitting employed for the 2D case. However, this splitting is actually slightly suboptimal. An improved splitting strategy involves the insertion of only three Steiner points (instead of four) and subdivides the triangular prism into ten tetrahedral elements (see Figure~\ref{prism_splitting}, C). This strategy is our foremost preference, as it produces a smaller number of elements relative to the first approach, while maintaining an identical pattern of splitting on the quadrilateral faces of the prism. Nevertheless, we make use of the more traditional splitting (the splitting into fourteen tetrahedra) in our subsequent numerical experiments due to its greater simplicity of implementation. The more optimal splitting (the splitting into ten tetrahedra) will be explored in future work.

For the sake of completeness, we introduce quantitative definitions for the aforementioned triangular prism splitting strategies on a reference triangular prism, denoted by $R^{\ast}$. We assume that $R^{\ast}$ has the following vertices
\begin{align*}
    R^{\ast} = \left[ r_1(0,0,0), \, r_2(1,0,0), \, r_3(0,1,0), \, r_4(0,0,1), \, r_5(1,0,1), \, r_6(0,1,1) \right].
\end{align*}
In addition, we introduce the following Steiner points at the centroid and on the quadrilateral faces of $R^{\ast}$,
\begin{align*}
    r_7 &= \frac{1}{4}(r_2+r_3+r_5+r_6), \quad r_8 =\frac{1}{4}(r_1 + r_2 +r_4+r_5), \\[1.5ex]
    r_9 &=\frac{1}{4} (r_1+r_3+r_4+r_6), \quad r_{10}= \frac{1}{6} \left(r_1 + r_2 + r_3 + r_4 + r_5 + r_6 \right).
\end{align*}
Based on the description above, we can define the following ten tetrahedra as part of subdivision strategy $\mathcal{C}$
\begin{align*}
    \mathcal{C}R^{\ast} = \Bigg\{&S_1(r_1, r_2, r_3, r_7), \, S_2(r_4, r_5, r_6, r_7), \, S_3(r_1, r_2, r_7, r_8), \\[1.5ex] &S_4(r_2, r_5, r_7, r_8), \, S_5(r_4, r_5, r_7, r_8), \, S_6(r_1, r_4, r_7, r_8), \\[1.5ex]
    &S_7(r_1, r_3, r_7, r_9), \, S_8(r_3, r_6, r_7, r_9), \, S_9(r_4, r_6, r_7, r_9), \, S_{10}(r_1, r_4, r_7, r_9) \Bigg\},
\end{align*}
where $\mathcal{C}$ yields the subdivision strategy illustrated in Figure~\ref{prism_splitting}, C. In addition, we can define the following fourteen tetrahedra as part of subdivision strategy~$\mathcal{E}$
\begin{align*}
    \mathcal{E}R^{\ast} = \Bigg\{&S_1(r_2, r_3, r_5, r_{10}), \, S_2(r_2, r_3, r_6, r_{10}), \, S_3(r_2, r_5, r_6, r_{10}), \\[1.5ex] &S_4(r_3, r_5, r_6, r_{10}), \, S_5(r_1, r_2, r_4, r_{10}), \, S_6(r_1, r_2, r_5, r_{10}), \\[1.5ex]
    &S_7(r_1, r_4, r_5, r_{10}), \, S_8(r_2, r_4, r_5, r_{10}), \, S_9(r_1, r_3, r_4, r_{10}), \\[1.5ex]
    &S_{10}(r_1, r_3, r_6, r_{10}), \, S_{11}(r_1, r_4, r_6, r_{10}), \, S_{12}(r_3, r_4, r_6, r_{10}), \\[1.5ex]
    &S_{13}(r_1, r_2, r_3, r_{10}), \, S_{14}(r_4, r_5, r_6, r_{10}) \Bigg\},
\end{align*}
where $\mathcal{E}$ yields the subdivision strategy illustrated in Figure~\ref{prism_splitting}, E.

Once the triangular prisms have been successfully subdivided into tetrahedra, then we recover a valid tetrahedral hypersurface mesh for the intermediate hypersurface. Thereafter, it remains for us to construct the hypersurface mesh on the terminating hyperplane. Towards this end, we collect the triangular elements and vertices associated with the terminating hyperplane (at $t = t_{n+1}$), then we send them off to a volume meshing program (such as Hang Si’s TetGen program~\cite{hang2015tetgen}). Once the terminating hypersurface mesh has been constructed, we synchronize its connectivity with the connectivity of the initial and intermediate hypersurface meshes. The resulting agglomeration of hypersurface meshes results in a hull of tetrahedral elements for the space-time slab from $t = t_n$ to $t = t_{n+1}$.

The process for generating the hypersurface mesh on a space-time slab is summarized below:

\begin{enumerate}
\item Extract the hypersurface mesh from the terminating hyperplane of the previous space-time slab.
\item Construct the hypersurface mesh for the initial hyperplane of the new space-time slab using the hypersurface mesh from step 1.
\item Extract the boundary triangular faces, edges, and vertices of the hypersurface mesh from step 2. Compute the vertex trajectories from $t = t_n$ to $t = t_{n+1}$. Project final point locations to the CAD surface.
\item Connect vertices on the terminating hyperplane to create triangular faces.
\item Connect triangles on the terminating hyperplane to triangles on the initial hyperplane to create triangular prisms.
\item Subdivide the triangular prisms into tetrahedra to generate a tetrahedral hypersurface mesh on the intermediate hypersurface.
\item Use the triangular faces on the terminating hyperplane to generate a tetrahedral hypersurface mesh on the terminating hyperplane.
\end{enumerate}

\begin{figure}[h!]
\centering
\includegraphics[width=13cm]{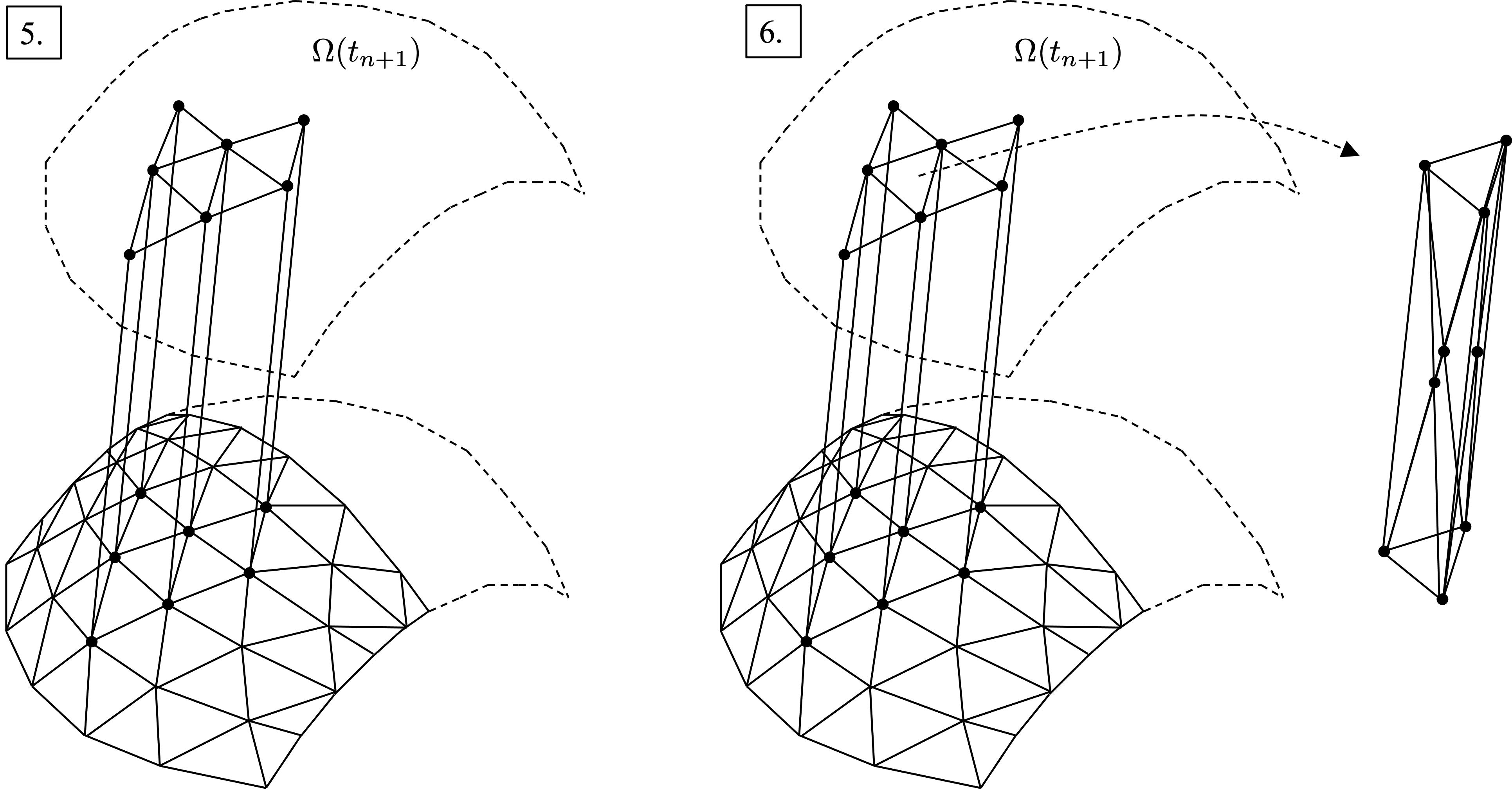}

\vspace{1cm}

\includegraphics[width=6.5cm]{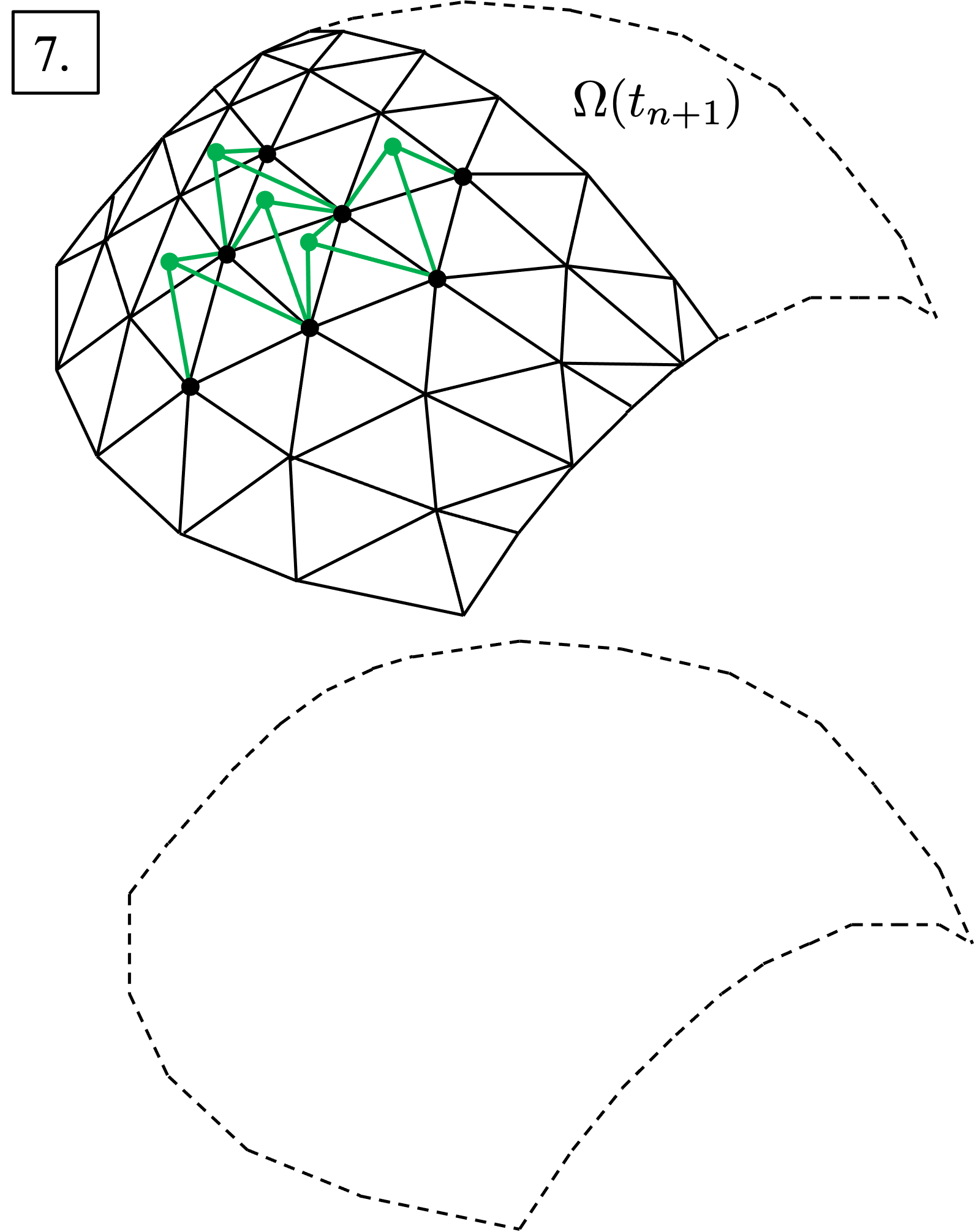}
\caption{An illustration of the process for generating a hypersurface mesh on a space-time slab in 3D+$t$. The numbered steps are explained in the text.}
\label{3D_surface_process}
\end{figure}

The hypersurface meshes from steps 2, 6, and 7 are combined to generate a hypersurface mesh for the entire space-time slab. The overall process is shown in Figure~\ref{3D_surface_process}.

\begin{figure}[h!]
\centering
\includegraphics[width=13cm]{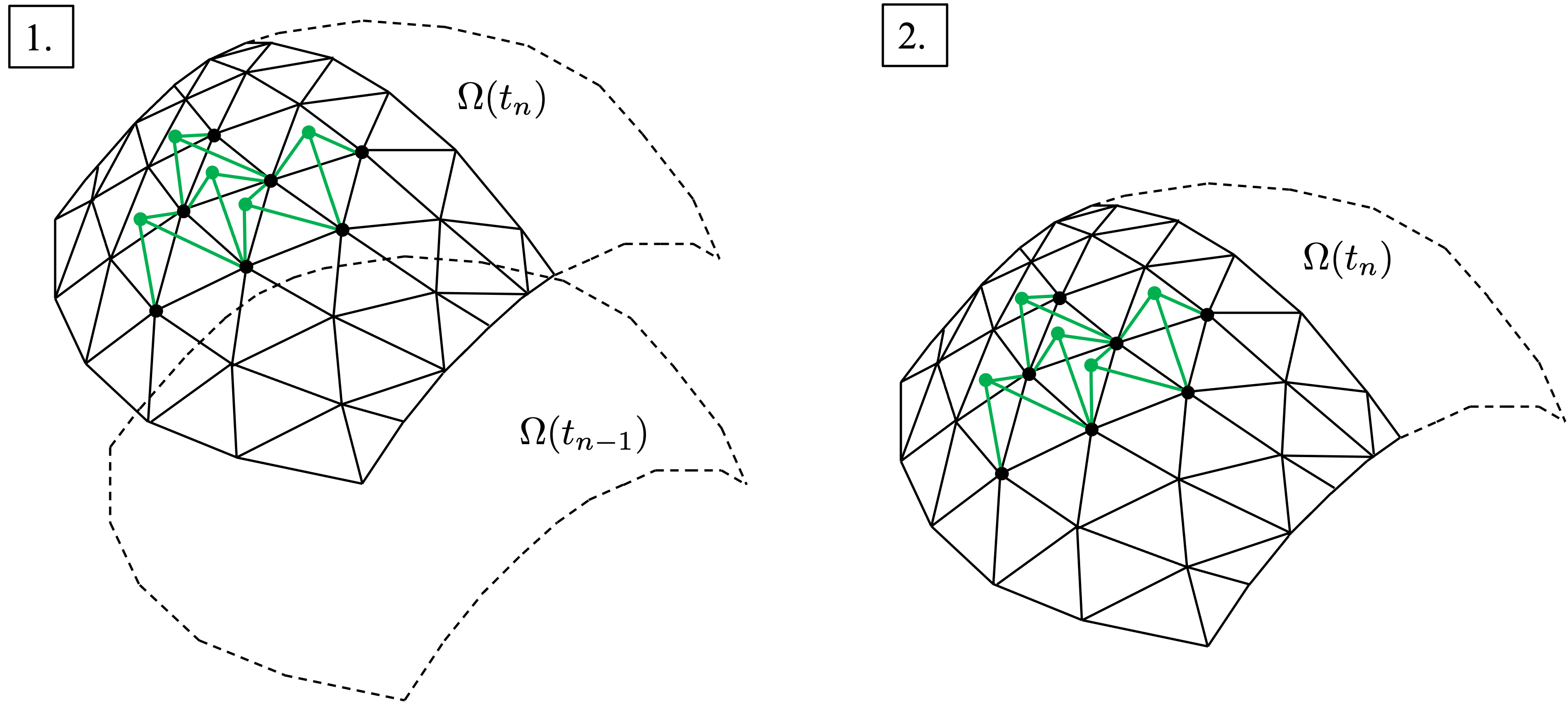}

\vspace{1cm}

\includegraphics[width=13cm]{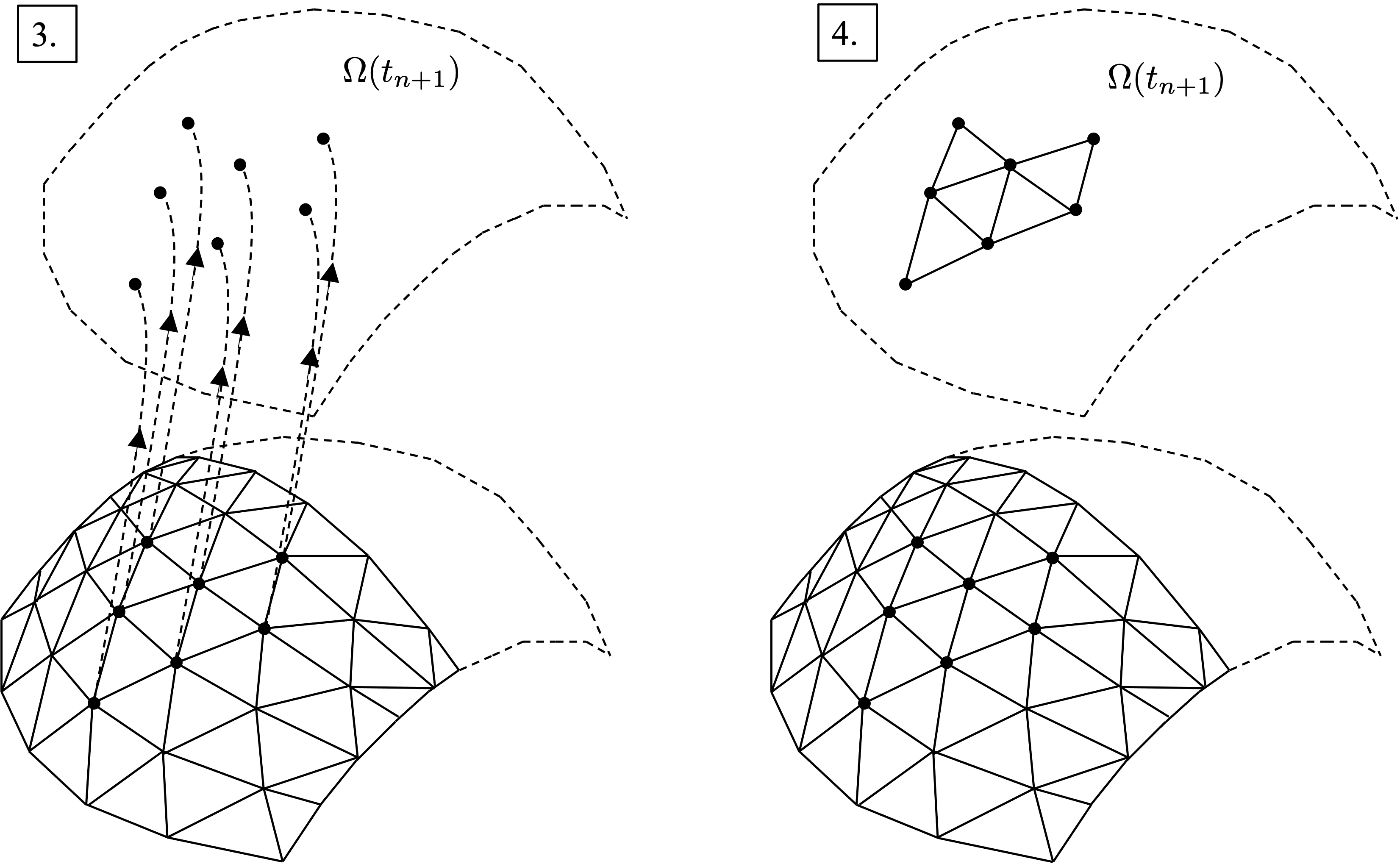}
\end{figure}

\pagebreak 
\clearpage

\section{Numerical Experiments}

In this section, we present numerical experiments using the surface meshing algorithm from the previous section. This algorithm was implemented as an extension of the JENRE$^{\text{\textregistered}}$  Multiphysics Framework used in earlier work for space-time finite element methods~\cite{Cor19_SCITECH}.

\subsection{Stationary Circle} \label{fixed_circle}

\begin{figure}[h!]
\centering
\includegraphics[width=8cm]{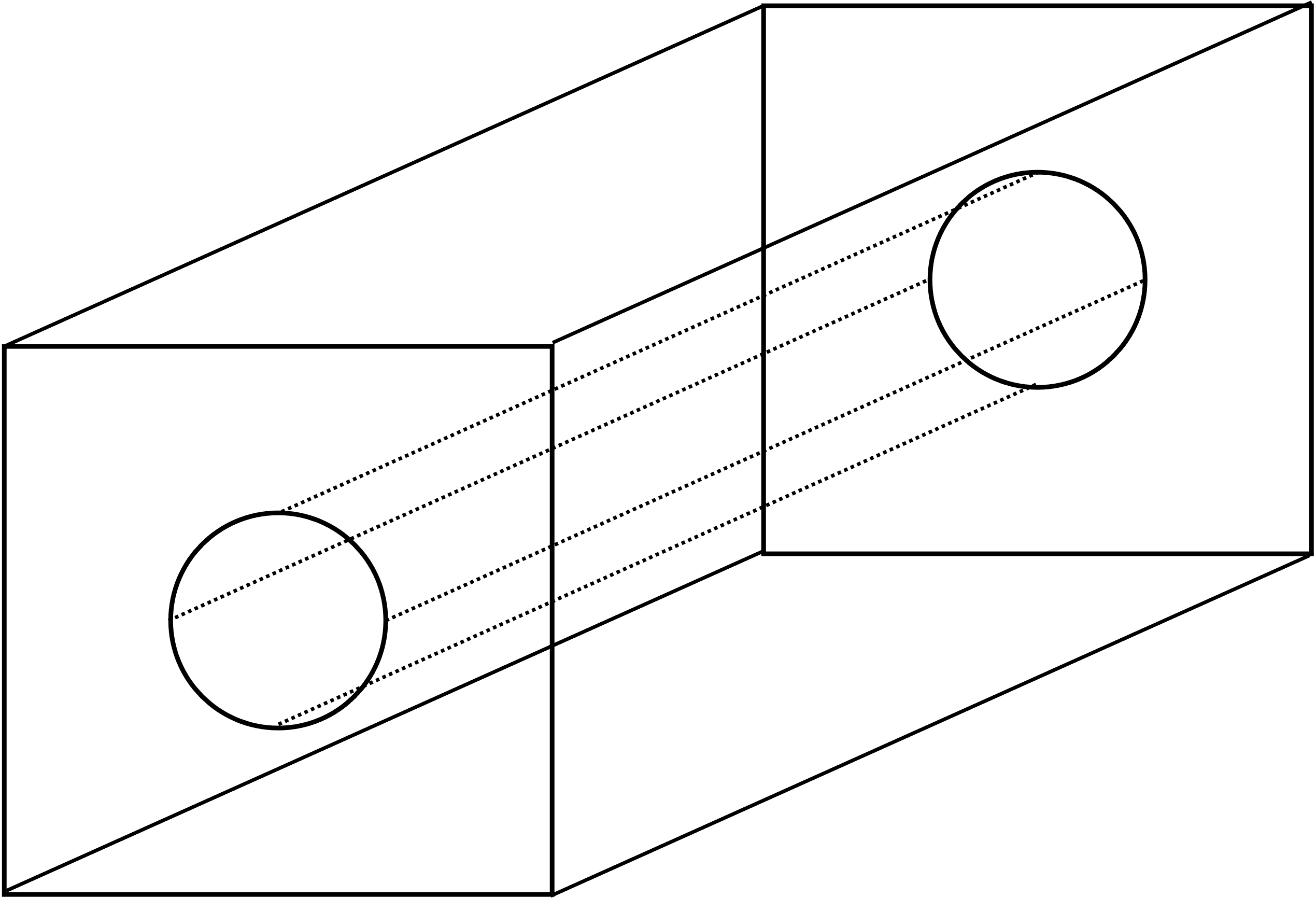}
\caption{An illustration of a stationary 2-sphere inside of a 2-cube. Under these circumstances, the space-time geometry consists of a cylinder embedded inside of a 3-cube. Note: this drawing is \emph{not} to scale.}
\label{simple_circle_fig}
\end{figure}

The spatial geometry for this test case consisted of a circle with constant radius $R = 1$, located inside of a square with constant edge length $L = 10$. The circle was positioned at the centroid of the square and was kept stationary during the time interval $t\in [0,1] = [t_0, t_f]$. The combination of the spatial domain and the temporal interval formed a space-time geometry consisting of a space-time cylinder inside a 3-cube. The geometry for this configuration is illustrated in Figure~\ref{simple_circle_fig}.
We note that this simple test case \emph{can} be treated by conventional mesh-extrusion methods. For example, a triangular surface mesh located at time $t_0 = 0$ can be extruded along the temporal direction to form a tetrahedral volume mesh that fills the space between the space-time cylinder and the surrounding cube. Furthermore, during this process, the boundary of the tetrahedral volume mesh automatically functions (or serves) as a triangular surface mesh which conforms to the space-time geometry. However, despite the simplicity of this test case and its ability to be successfully treated with other methods, it nonetheless serves as a useful `sanity-check' in order to ensure that our surface meshing algorithm is working as expected. With this justification in mind, we proceeded by constructing a preliminary triangular surface mesh for the spatial geometry at $t_0 =0$. Thereafter, we constructed a family of surface meshes for the entire space-time domain using the techniques from the previous section. For each of these surface meshes, the characteristic element size near the circle, $h_{\text{circle}}$, and the characteristic element size near the square boundary, $h_{\text{square}}$, were specified on the initial surface mesh at $t_0 = 0$. In addition, the mesh spacing of the domain along the temporal direction, $h_{\text{time}}$, was specified. Next, a total of nine surface meshes for the space-time slab were generated, each with a greater number of elements than the previous mesh in the sequence. Note that the surface meshes that appeared later in the sequence were larger than the earlier ones because they had progressively smaller values of $h_{\text{circle}}$, $h_{\text{square}}$, and $h_{\text{time}}$. These spacing parameters were usually decreased by a factor of between 1.25 and 2.0 between successive meshes. The essential properties of the resulting surface meshes are summarized in Table~\ref{circle_stat_table}. 
\begin{table}[h!]
\begin{center}
\begin{tabular}{| c |r|r| }
\hline
Mesh & Elements & Vertices \\
 \hline
1& 2,714 &	1,357 \\
2& 5,308 &	2,654 \\
3& 10,022 &	5,011 \\
4& 20,142 &	10,071 \\
5& 38,792 &	19,396 \\
6& 76,524 &	38,262 \\
7& 154,570 & 77,285 \\
8 & 305,602 & 152,801 \\
9 & 613,418 & 306,709 \\
\hline
\end{tabular}
\caption{The number of triangular elements and vertices for a sequence of surface meshes for the stationary circle test case.} \label{circle_stat_table}
\end{center}
\end{table}

The validity of each surface mesh was assessed by comparing its approximate surface area to the exact, analytically-determined surface area of the space-time geometry. The approximate surface area for each mesh was calculated by summing the areas of all triangles in each mesh. In particular, the area of each individual triangle was calculated using Heron's formula,
\begin{align*}
    A_{\text{approx}} &= \sum_{T_k} A_k,
\end{align*}
where $T_k$ is a generic triangle in a given surface mesh and
\begin{align*}   
    A_k &= \sqrt{s_k(s_k - a_k)(s_k - b_k)(s_k - c_k)}, \qquad s_k = \frac{1}{2}\left(a_k + b_k + c_k\right),
\end{align*}
where $a_k$, $b_k$, and $c_k$ are the edge lengths of the $k$-th triangle. 

\begin{figure}[h!]
\centering
\includegraphics[width=8cm]{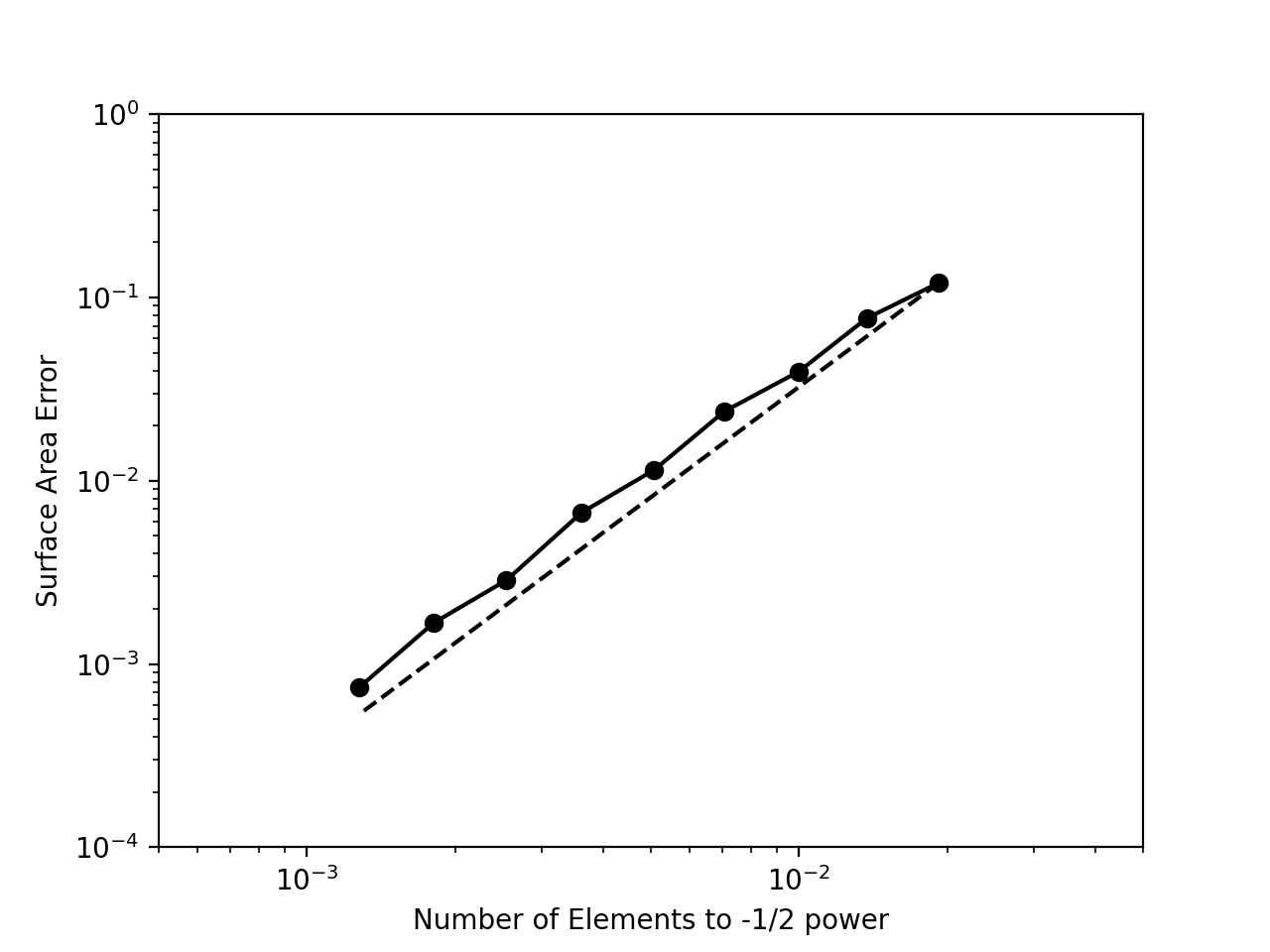}
\caption{Each point on the plot above represents the error between the area of a surface mesh and the exact surface area for the stationary circle test case. The errors are plotted against the characteristic mesh spacing for a sequence of increasingly refined surface meshes. In addition, a dashed line associated with 2nd-order convergence is plotted for reference.}
\label{stat_circle_error_fig}
\end{figure}

The exact surface area of the space-time geometry was calculated by the following formula
\begin{align*}
    A_{\text{exact}} = 2 \left( L^2 - \pi R^2\right) + \left(4 L + 2 \pi R  \right)\left(t_f- t_0 \right).
\end{align*}
In a natural fashion, the error was calculated as follows
\begin{align*}
    A_{\text{error}} = \left| A_{\text{exact}} - A_{\text{approx}} \right|.
\end{align*}
Figure~\ref{stat_circle_error_fig} shows a plot of the error in the surface area versus the number of elements to the -1/2 power. Here, we can see that the error decays at a rate of 2nd-order as the mesh resolution increases. This rate of convergence agrees well with our expectations, as straight-sided triangular elements are expected to generate 2nd-order convergence rates for most finite element applications. 

\subsection{Expanding Circle}

\begin{figure}[h!]
\centering
\includegraphics[width=8cm]{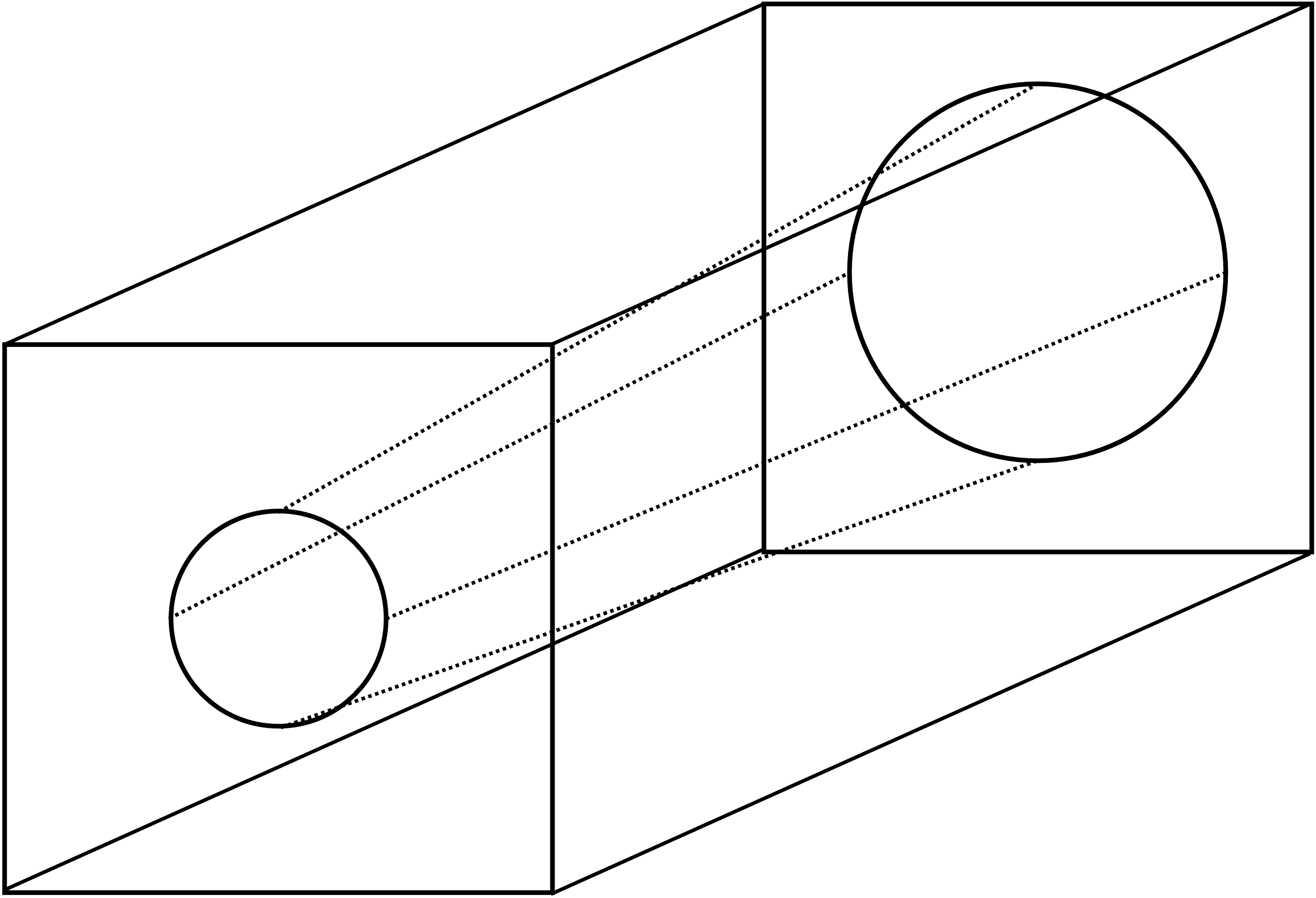}
\caption{An illustration of an expanding 2-sphere inside of a 2-cube. Under these circumstances, the space-time geometry consists of a conical frustum embedded inside of a 3-cube. Note: this drawing is \emph{not} to scale.}
\label{expand_circle_fig}
\end{figure}

\begin{figure}[h!]
\centering
\includegraphics[width=8cm]{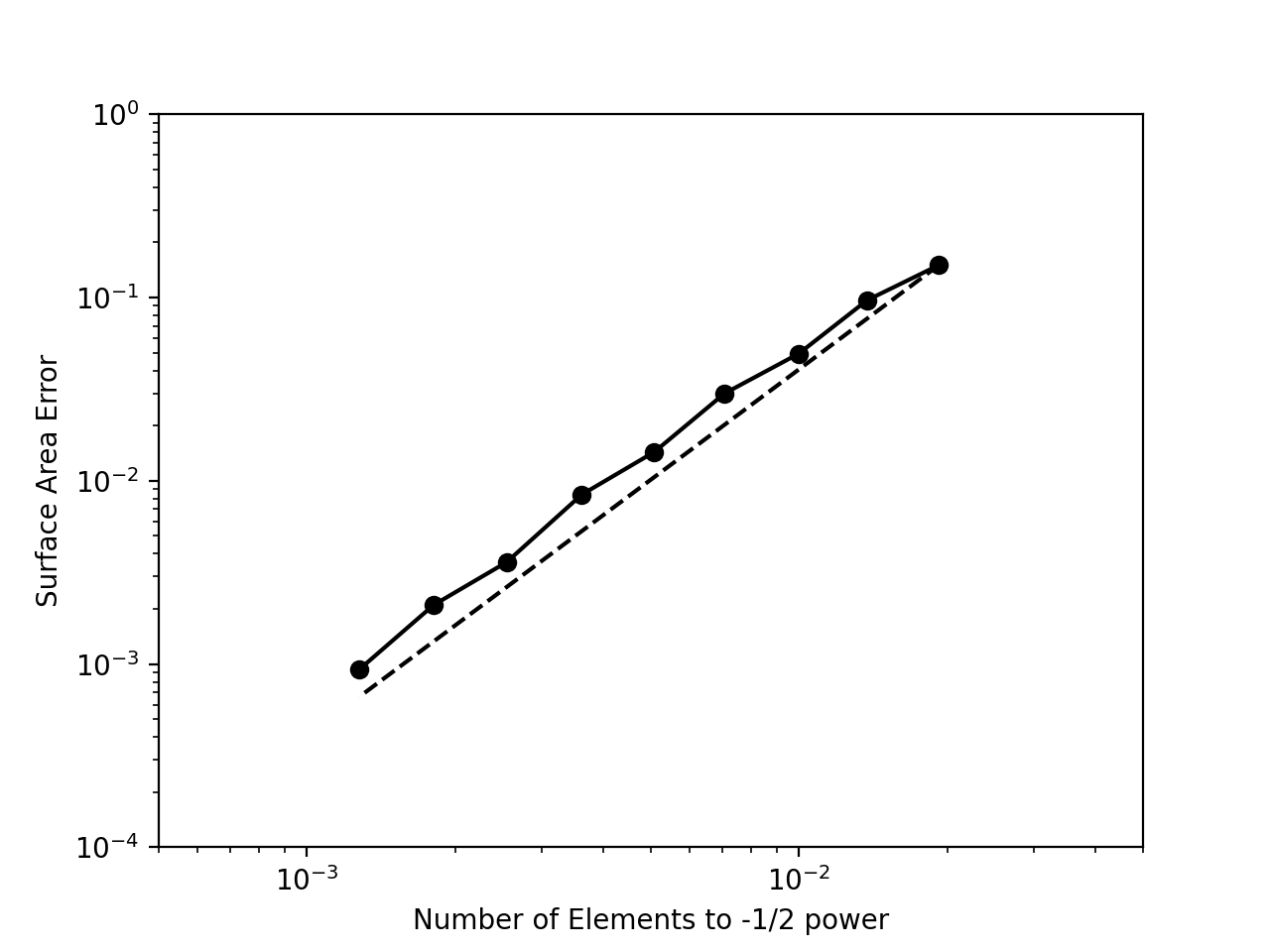}
\caption{Each point on the plot above represents the error between the area of a surface mesh and the exact surface area for the expanding circle test case. The errors are plotted against the characteristic mesh spacing for a sequence of increasingly refined surface meshes.  In addition, a dashed line associated with 2nd-order convergence is plotted for reference.}
\label{expand_circle_error_fig}
\end{figure}

In this test case, the circle from the previous case was allowed to expand. In particular, the radius of the circle was calculated based on the following function
\begin{align}
    R(t) = m t + R_0,
    \label{radius_formula}
\end{align}
where $R_0$ is the initial radius of the circle, and $m$ is the radial expansion speed of the circle. In this case, we elected to set $R_0 = 1$ and $m = 0.25$, and we allowed the circle to expand during the time interval $[0,1]$. The final radius of the circle was $R_f = 1.25$. The space-time geometry for this case is a conical frustum with initial radius $R_0$ and final radius $R_f$ inside of a 3-cube with edge length $L = 10$. In a natural fashion, the axis of revolution for the conical frustum is aligned with the temporal axis. Figure~\ref{expand_circle_fig} shows an illustration of this geometric configuration. 
We created a family of nine surface meshes for the chosen space-time geometry, using the meshing parameters and techniques described in Section~\ref{fixed_circle}. The properties of the surface meshes for the expanding circle are summarized in Table~\ref{circle_expand_table}. 
\begin{table}[h!]
\begin{center}
\begin{tabular}{| c |r|r| }
\hline
Mesh & Elements & Vertices \\
\hline
1& 2,696 & 1,348 \\
2& 5,288 & 2,644 \\
3& 10,016 & 5,008 \\
4& 20,086 & 10,043 \\
5& 38,734 &	19,367 \\
6& 76,406 &	38,203 \\
7& 154,362 & 77,181 \\
8& 305,194 & 152,597 \\
9& 612,452 & 306,226 \\
\hline
\end{tabular}
\caption{The number of triangular elements and vertices for a sequence of surface meshes for the expanding circle test case.} \label{circle_expand_table}
\end{center}
\end{table}
We assessed the validity of the surface meshes by calculating the area of each mesh, and comparing it with the following analytically-determined exact area for the space-time slab
\begin{align*}
    A_{\text{exact}} &= 2L^2 - \pi (R_f^2 + R_0^2) + 4L\left(t_f - t_0\right) \\
    &+ \pi (R_f + R_0) \sqrt{(R_f - R_0)^2 + (t_f - t_0)^2}.
\end{align*}
Figure~\ref{expand_circle_error_fig} shows a plot of the surface area error versus the approximate mesh spacing. As expected, the error decreases at a rate of 2nd-order with increasing mesh resolution.

\subsection{Stationary Sphere}

\begin{figure}[h!]
\centering
\includegraphics[width=11cm]{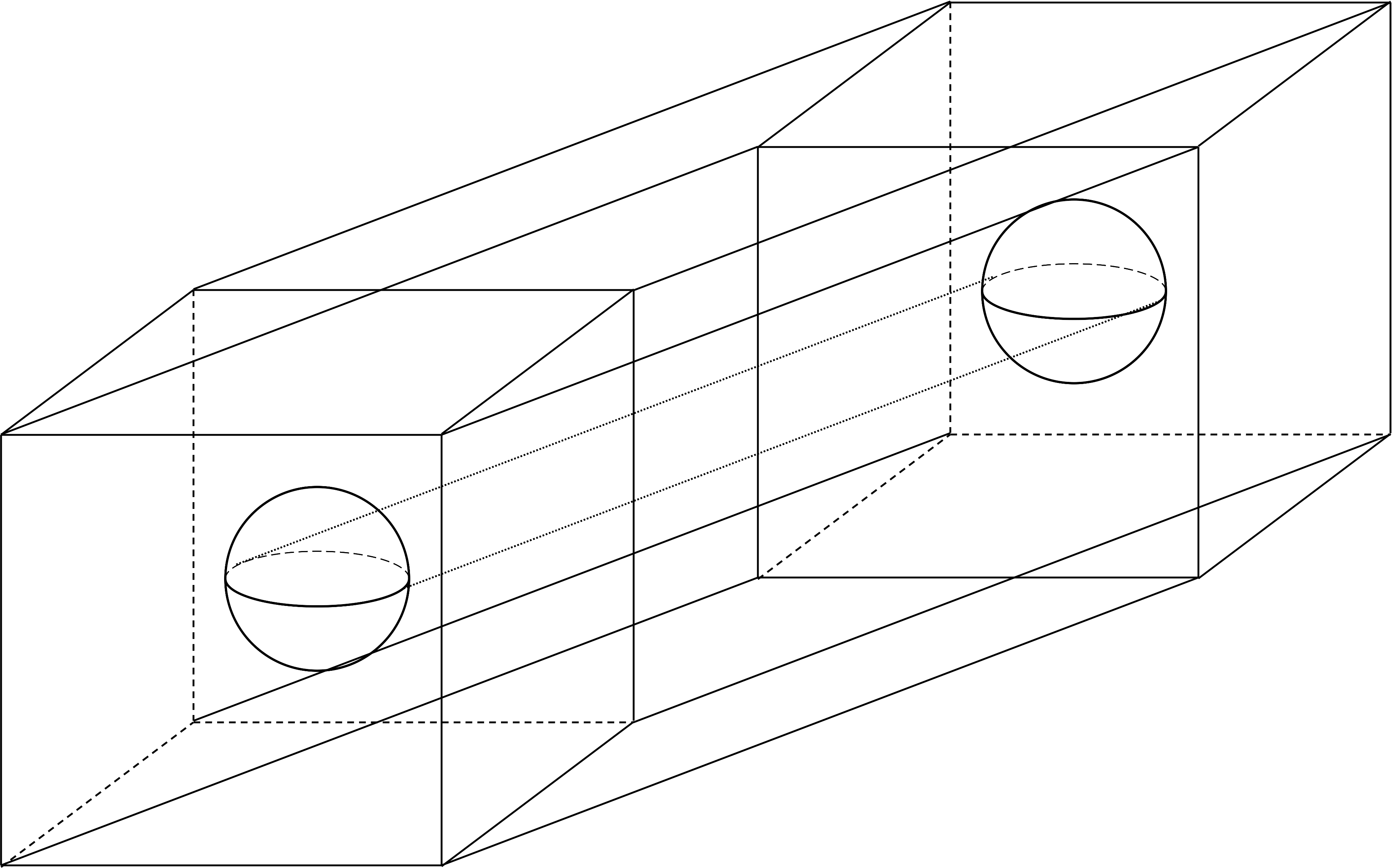}
\caption{An illustration of a stationary 3-sphere inside of a 3-cube. Under these circumstances, the space-time geometry consists of a hypercylinder embedded inside of a tesseract. Note: this drawing is \emph{not} to scale.}
\label{simple_sphere_fig}
\end{figure}

The geometry for this experiment consisted of a stationary 3-sphere with radius $R = 1$ inside of a 3-cube with edge length $L = 10$. The sphere was located at the centroid of the cube. In addition, the surface of the sphere was kept static without any changes in size. The associated space-time geometry consists of a hypercylinder embedded inside of a tesseract (4-cube). This geometry is shown in Figure~\ref{simple_sphere_fig}.
The region between the surface of the sphere and the walls of the cube was filled with an unstructured mesh of tetrahedral elements at time $t_0 = 0$. This mesh served as a hypersurface mesh for the initial hyperplane. 
With this as a starting point, an entire family of hypersurface meshes was formed for the space-time slab using the construction techniques described in Section~\ref{surface_mesh_approach}. In order to create a well-behaved family of meshes, we specified the mesh spacings on the surface of the sphere, $h_{\text{sphere}}$, on the surface of the cube walls, $h_{\text{cube}}$, and along the temporal direction, $h_{\text{time}}$. 
The mesh properties are summarized in Table~\ref{sphere_stat_table}.
\begin{table}[h!]
\begin{center}
\begin{tabular}{| c |r|r| }
\hline
Mesh & Elements & Vertices \\
\hline
1& 144,431 &	30,813 \\
2& 368,587 &	78,637 \\
3& 958,500 &	204,440 \\
4& 2,630,598 & 561,013 \\
5& 7,127,597 & 1,519,319 \\
6& 19,165,615 & 4,087,387 \\
7& 56,118,477 & 11,961,783 \\
8& 149,470,428 & 31,879,852 \\
\hline
\end{tabular}
\caption{The number of tetrahedral elements and vertices for a sequence of hypersurface meshes for the stationary sphere test case.} \label{sphere_stat_table}
\end{center}
\end{table}

We compared the volume of each hypersurface mesh to the exact, analytically-determined volume of the hypersurface for the space-time slab. The volume of each hypersurface mesh was calculated by adding up the individual volumes of all tetrahedral elements in each mesh as follows
\begin{align*}
    V_{\text{approx}} = \sum_{T_k} V_k,
\end{align*}
where
\begin{align*}
    V_k = \sqrt{\frac{\text{det}(\Theta)}{288}}, \qquad \Theta = \begin{bmatrix} 0& 1 & 1& 1& 1 \\
    1 & 0 & d_{ab}^{2} & d_{ac}^{2} & d_{ae}^{2} \\
    1 & d_{ab}^{2} & 0 & d_{bc}^{2} & d_{be}^{2} \\
    1 & d_{ac}^{2} & d_{bc}^{2} & 0 & d_{ce}^{2} \\
    1 & d_{ae}^{2} & d_{be}^{2} & d_{ce}^{2} & 0\end{bmatrix},
\end{align*}
and where the ``$d$" quantities above are the pairwise distances between the vertices $\bm{a}$, $\bm{b}$, $\bm{c}$, and $\bm{e}$ of the $k$-th tetrahedron, $T_k$. 

\begin{figure}[h!]
\centering
\includegraphics[width=8cm]{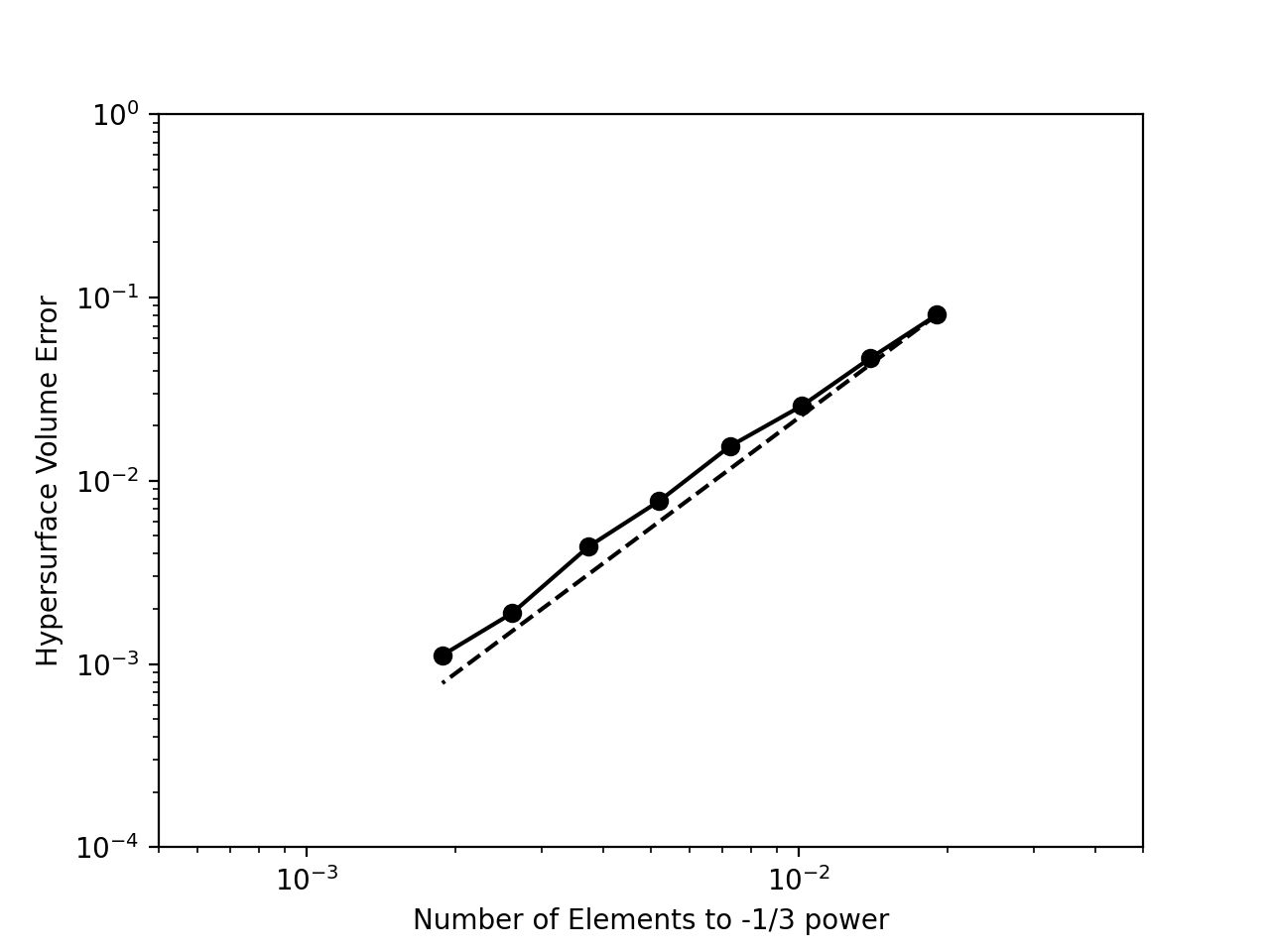}
\caption{Each point on the plot above represents the error between the volume of a hypersurface mesh and the exact volume for the stationary sphere test case. The errors are plotted against the characteristic mesh spacing for a sequence of increasingly refined hypersurface meshes. In addition, a dashed line associated with 2nd-order convergence is plotted for reference.}
\label{stat_sphere_error_fig}
\end{figure}

The analytically-determined, exact volume was computed as follows
\begin{align*}
    V_{\text{exact}} = 2\left( L^{3} - \frac{4}{3} \pi R^3 \right) + \left(6L^{2} + 4\pi R^{2}\right)(t_f - t_0).
\end{align*}
The error in the hypersurface volume was then obtained as follows
\begin{align*}
    V_{\text{error}} = \left| V_{\text{exact}} - V_{\text{approx}} \right|.
\end{align*}
Figure~\ref{stat_sphere_error_fig} shows the volumetric error for each hypersurface mesh plotted versus the approximate mesh spacing. Here, the mesh spacing was estimated by raising the total number of elements in each mesh to the -1/3 power. As expected, the error appears to consistently decrease with a rate of approximately 2nd order.

\subsection{Expanding Sphere}

\begin{figure}[h!]
\centering
\includegraphics[width=11cm]{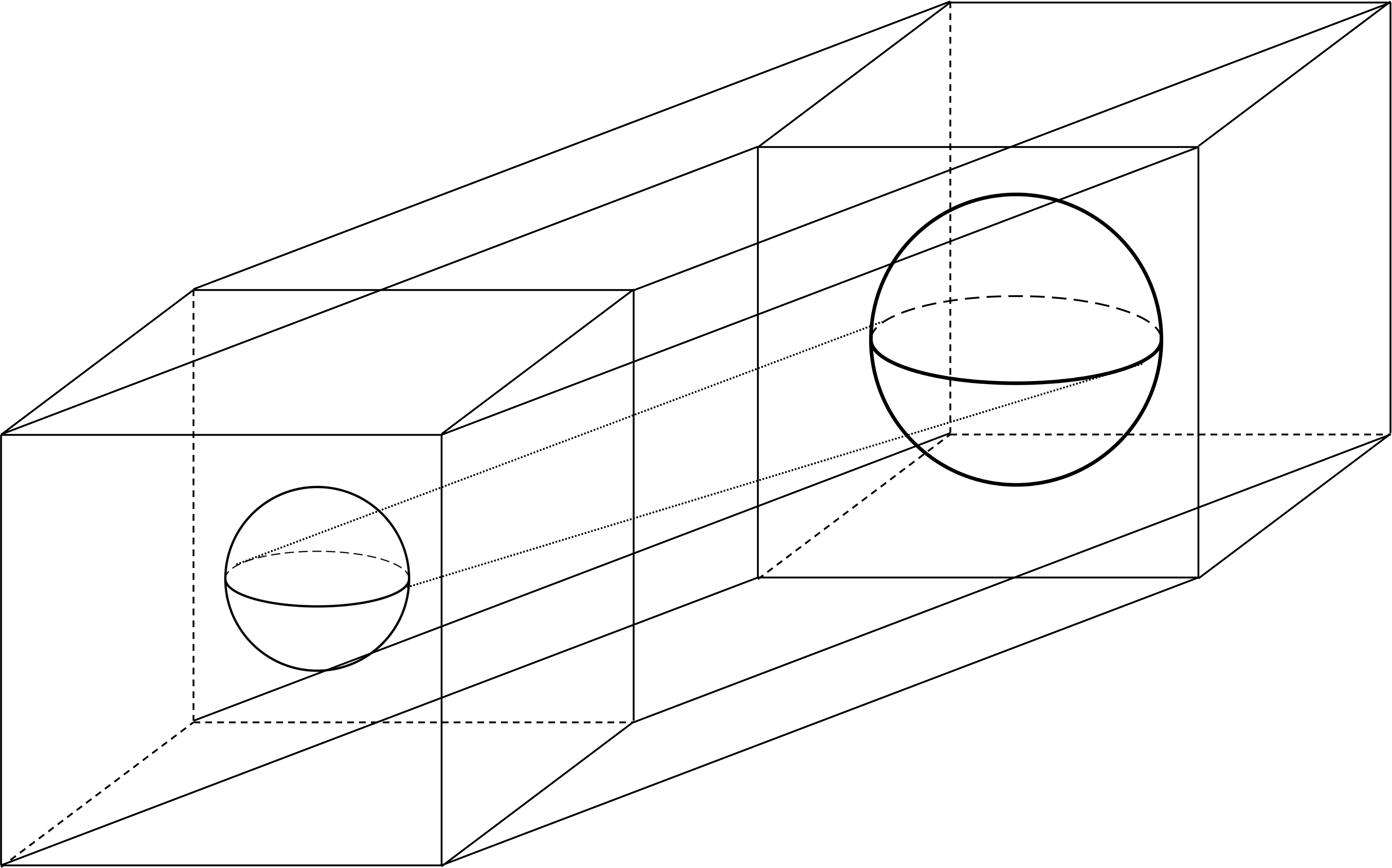}
\caption{An illustration of an expanding 3-sphere inside of a 3-cube. Under these circumstances, the space-time geometry consists of a hyper-conical frustum embedded inside of a tesseract. Note: this drawing is \emph{not} to scale.}
\label{expand_sphere_fig}
\end{figure}

\begin{figure}[h!]
\centering
\includegraphics[width=8cm]{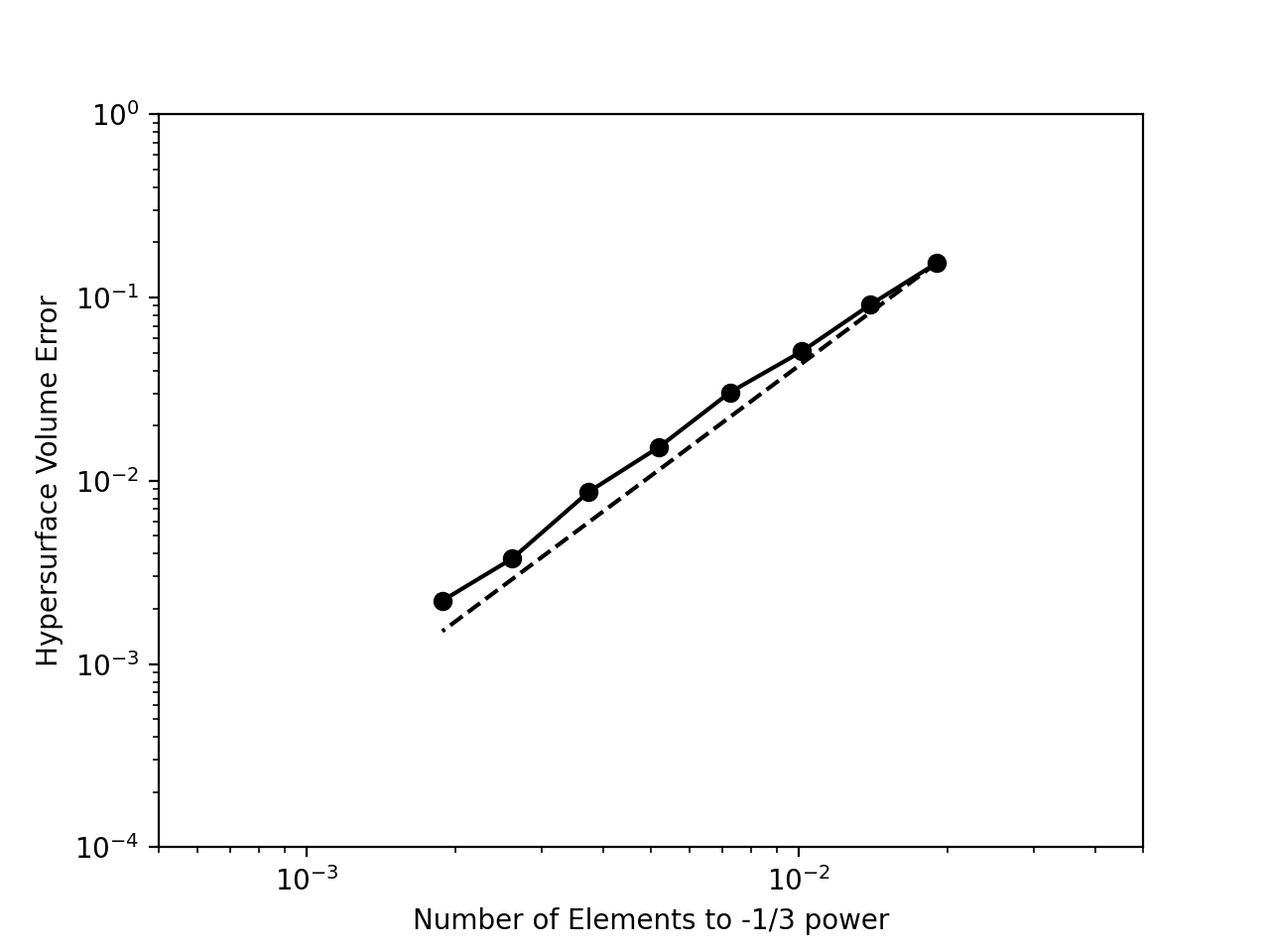}
\caption{Each point on the plot above represents the error between the volume of a hypersurface mesh and the exact volume for the expanding sphere test case. The errors are plotted against the characteristic mesh spacing for a sequence of increasingly refined hypersurface meshes. In addition, a dashed line associated with 2nd-order convergence is plotted for reference.}
\label{expand_sphere_error_fig}
\end{figure}

For this experiment, we used the stationary sphere geometry from the previous section. However, in this case, the radius of the sphere was allowed to increase in time in accordance with Eq.~\eqref{radius_formula}, during the time interval $[0,1]$. Here, we let $R_0 = 1$ and $m = 0.25$. During the time interval in question, the sphere expanded to a final radius of $R_f = 1.25$. The associated space-time geometry consisted of a hyper-conical frustum embedded inside of a tesseract. This geometry is shown in Figure~\ref{expand_sphere_fig}.
For this geometry, we created a family of hypersurface meshes using the procedure described in the previous section. The mesh properties are summarized in Table~\ref{sphere_expand_table}.
\begin{table}[h!]
\begin{center}
\begin{tabular}{| c |r|r| }
\hline
Mesh & Elements & Vertices \\
\hline
1& 144,345 &	30,801 \\
2& 368,477 &	78,620 \\
3& 958,364 &	204,417 \\
4& 2,630,263 & 560,955 \\
5& 7,126,629 & 1,519,159 \\
6& 19,163,704 & 4,087,112 \\
7& 56,114,273 & 11,961,083 \\
8& 149,461,495 & 31,878,349 \\
\hline
\end{tabular}
\caption{The number of tetrahedral elements and vertices for a sequence of hypersurface meshes for the expanding sphere test case.} \label{sphere_expand_table}
\end{center}
\end{table}
The total volume of each mesh was compared with the exact volume of the slab's hypersurface, which was computed as follows
\begin{align*}
    V_{\text{exact}} &= 2 L^{3} - \frac{4}{3} \pi \left(R_f^{3} + R_{0}^{3} \right) + 6 L^{2} \left(t_f - t_0\right) \\
    &+ \frac{4}{3} \pi \left( \frac{R_f^3 - R_0^3}{R_f - R_0} \right) \sqrt{(R_f - R_0)^{2} + (t_f - t_0)^{2}}.
\end{align*}
Figure~\ref{expand_sphere_error_fig} shows a plot of the volumetric error versus the approximate mesh spacing. As expected, the error in the approximation deceases with increasing mesh resolution, and the rate of decrease is approximately 2nd order.

\subsection{Rotating Ellipsoid}

\begin{figure}[h!]
\centering
\includegraphics[width=6cm,trim={2cm 0 2cm 0},clip]{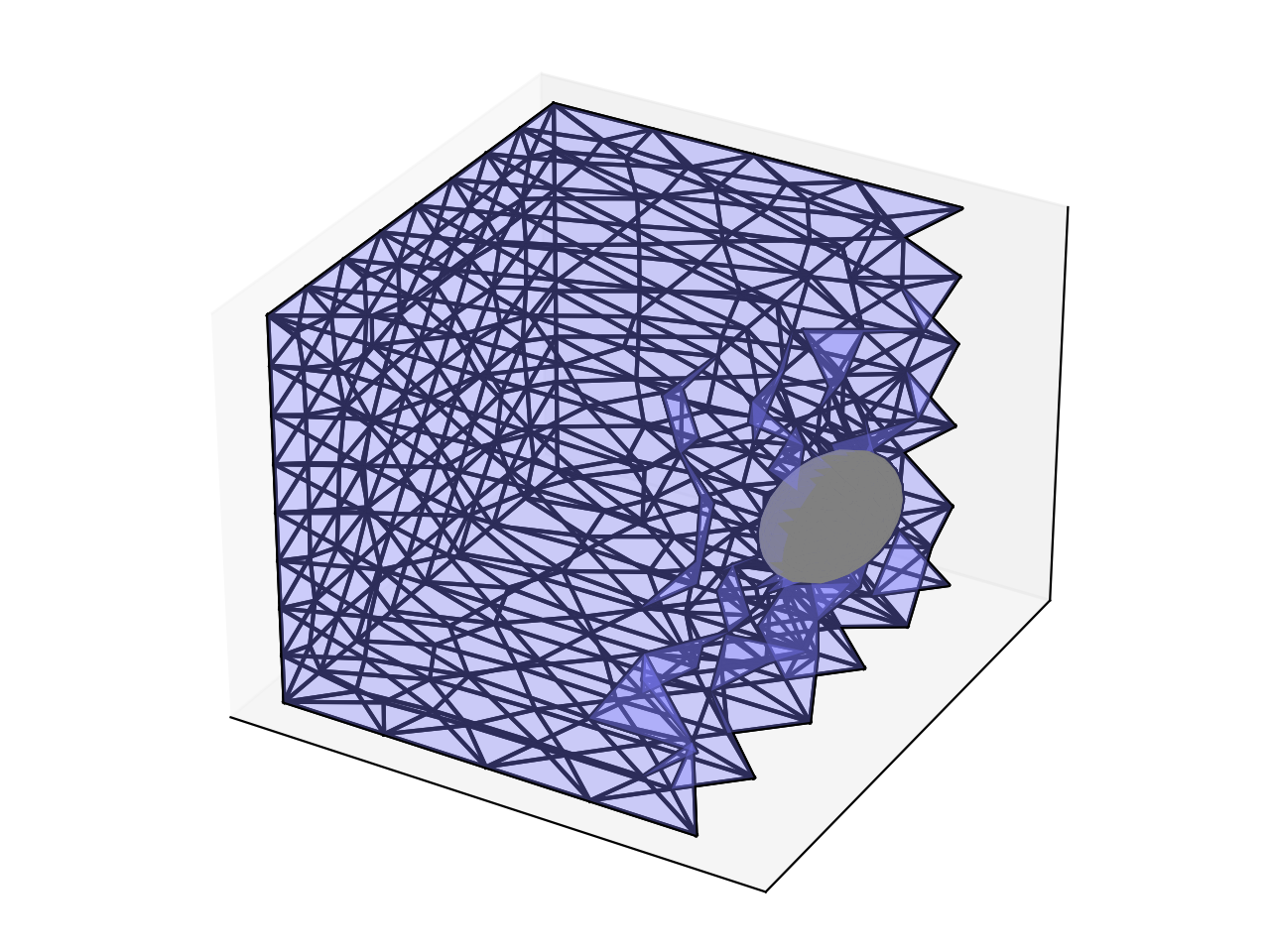} 
\includegraphics[width=6cm,trim={2cm 0 2cm 0},clip]{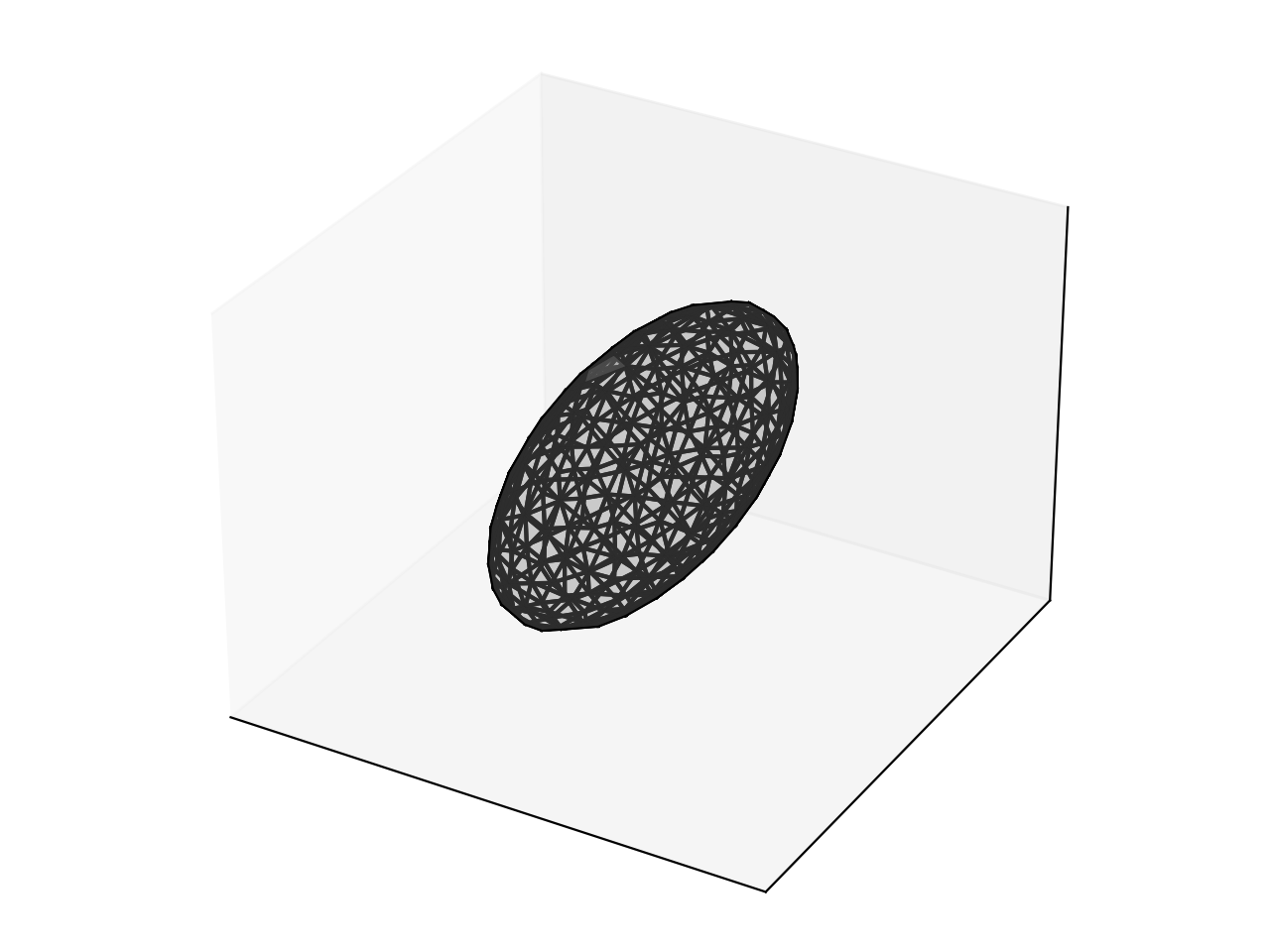}
\includegraphics[width=6cm,trim={2cm 0 2cm 0},clip]{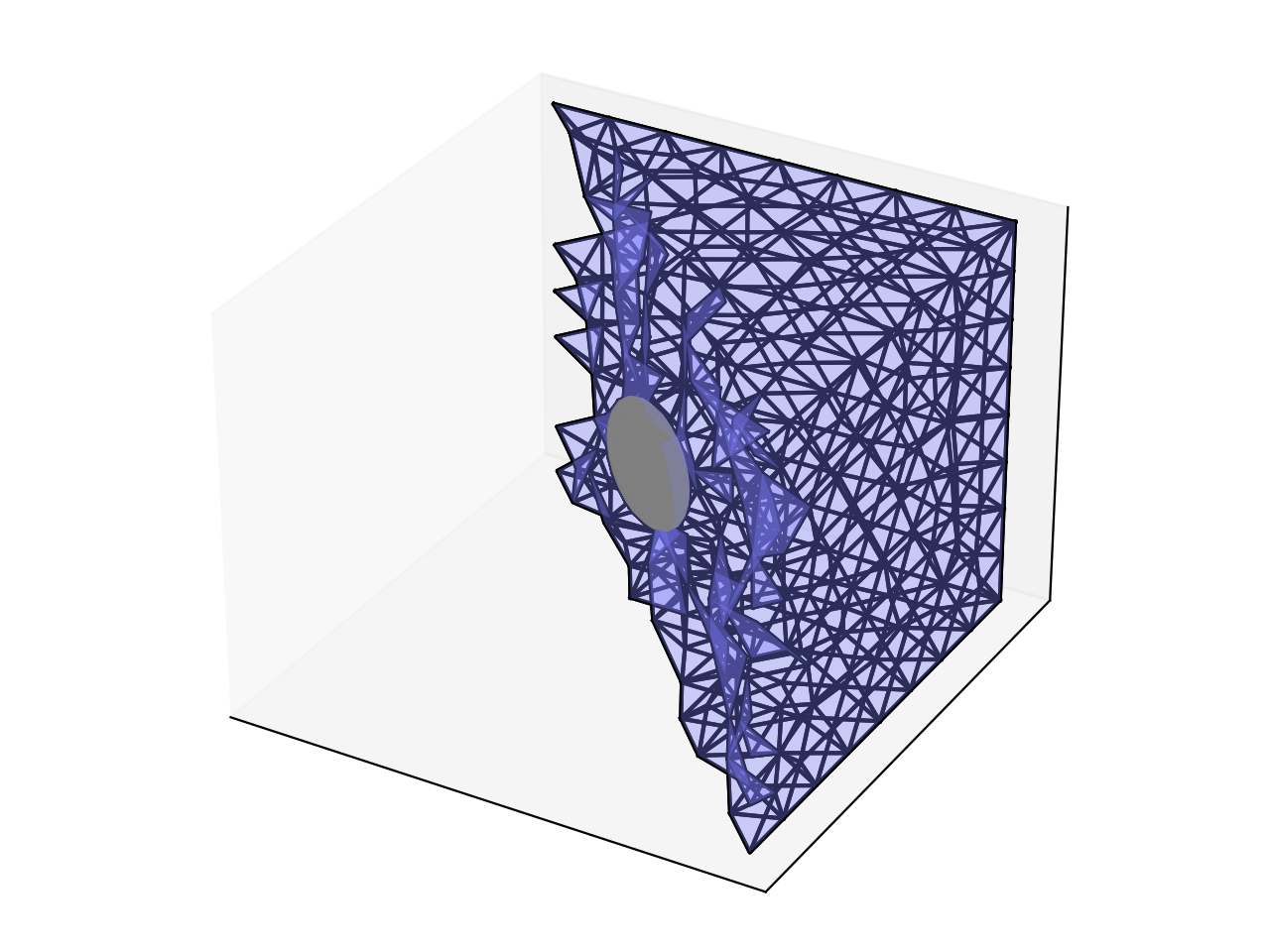} 
\includegraphics[width=6cm,trim={2cm 0 2cm 0},clip]{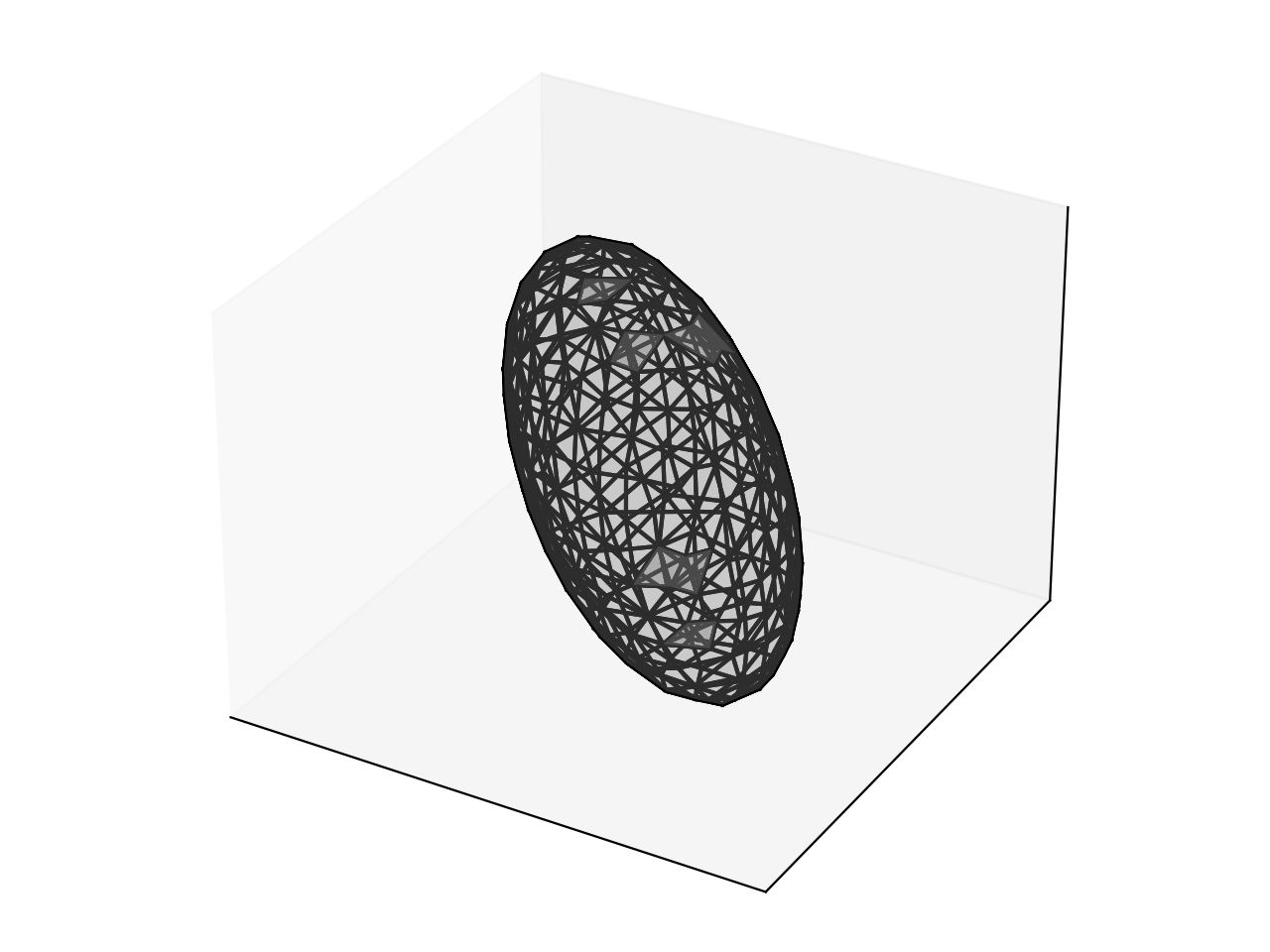}
\includegraphics[width=6cm,trim={2cm 0 2cm 0},clip]{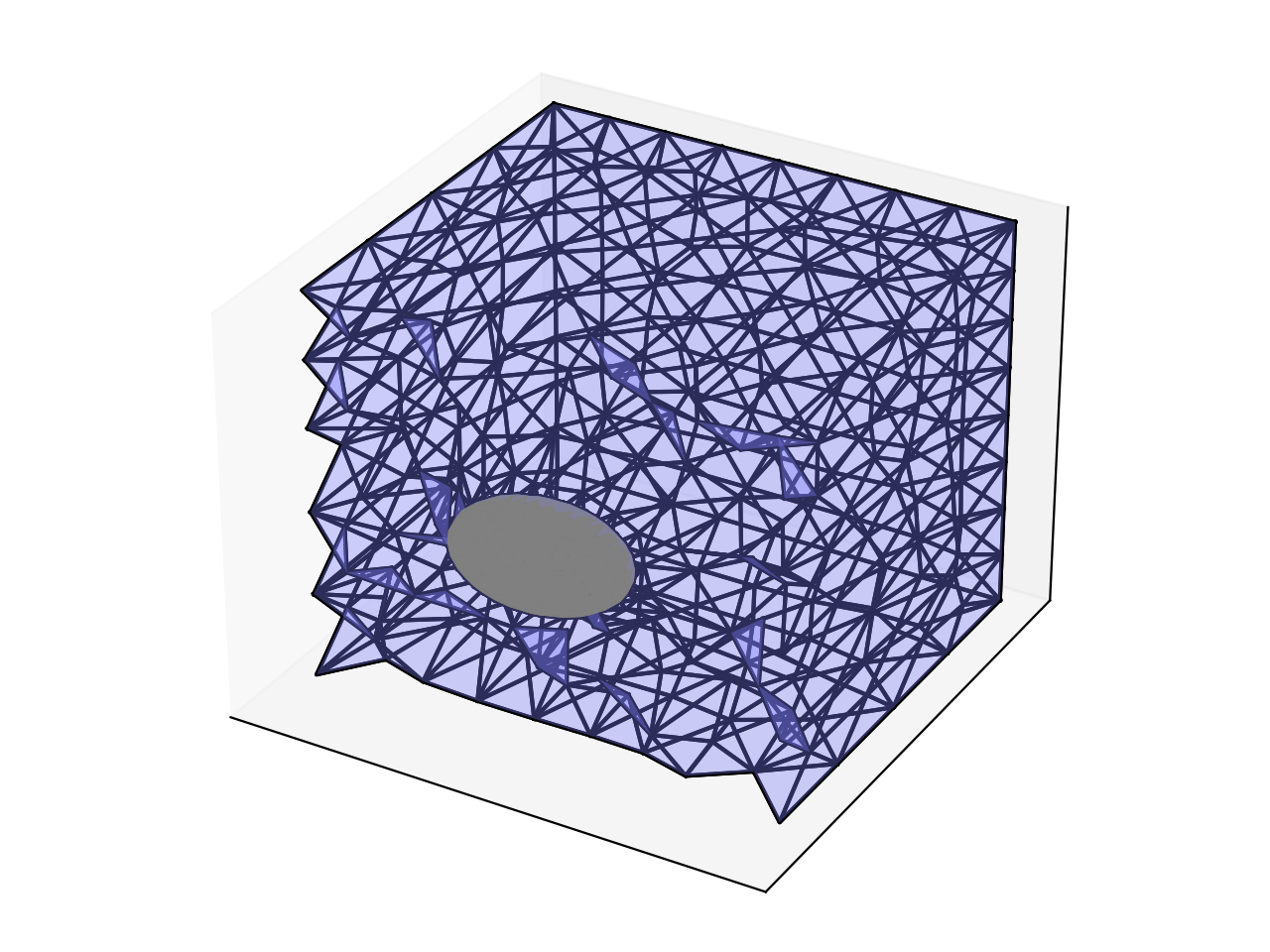} 
\includegraphics[width=6cm,trim={2cm 0 2cm 0},clip]{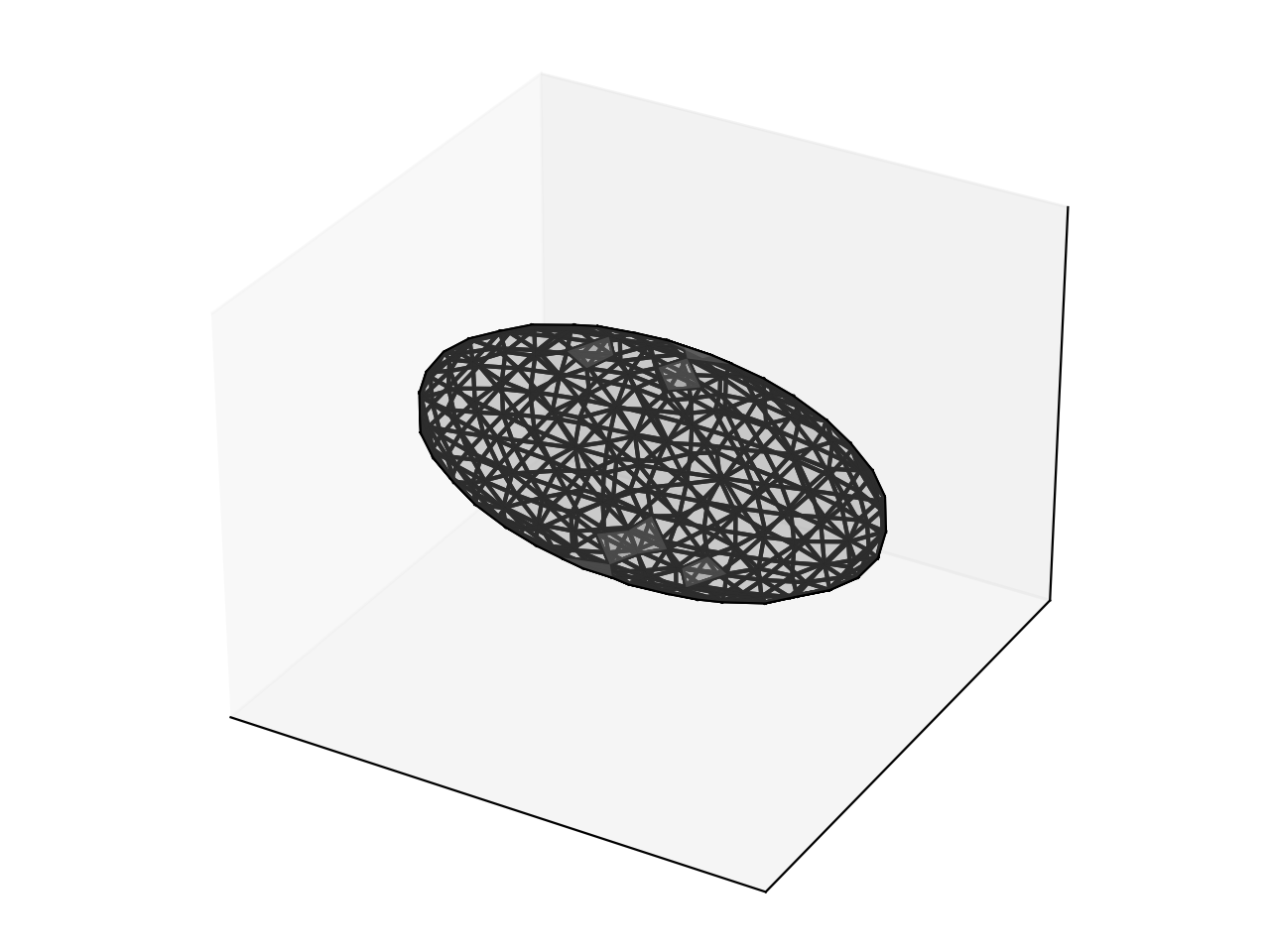}
\caption{Snapshots of a rotating ellipsoid at $t = 0$, $t = 0.5$, and $t = 1$, (top to bottom). On the left, a cross section of the coarsest hypersurface mesh is shown alongside the ellipsoid CAD. On the right, a zoomed-in view of the surface mesh on the ellipsoid is shown.}
\label{rotating_ellipsoid_snapshots}
\end{figure}

The geometry for this test case consisted of a single ellipsoid with semi-axes of $a = 1$ in the $x$-direction, $b = 3$ in the $y$-direction, and $c = 2$ in the $z$-direction. The ellipsoid was located at the center of a 3-cube with edge length $L = 16$. In this 3-cube, the ellipsoid rotated around the $z$-axis with a constant speed of $\omega = \frac{\pi}{2}$ rads/s, during the time interval $[0,1]$. The resulting space-time geometry consisted of an `ellipsoidal hyper-helix' contained inside of a tesseract. With this geometric configuration in mind, we generated a family of hypersurface meshes using the procedure described in the previous section. The hypersurface meshes were parameterized by the following quantities: $h_{\text{ellipsoid}}$ was used to specify the mesh spacing near the ellipsoid surface, $h_{\text{cube}}$ was used to specify the spacing near the cube walls, and $h_{\text{time}}$ was used to specify the spacing along the temporal direction.  The properties of the resulting hypersurface meshes are summarized in Table~\ref{ellipsoid_table}. In addition, Figure~\ref{rotating_ellipsoid_snapshots} shows some representative snapshots of the coarsest (lowest-resolution) hypersurface mesh.
\begin{table}[h!]
\begin{center}
\begin{tabular}{| c |r|r| }
\hline
Mesh & Elements & Vertices \\
\hline
1 & 355,798 & 75,655 \\
2 & 944,622 & 200,867 \\
3 & 2,644,079 & 562,277 \\
4 & 7,236,458 & 1,538,687 \\
5 & 20,536,468 & 4,365,033 \\
6 & 54,872,975 & 11,676,192 \\
7 & 162,044,418 & 34,449,255 \\
\hline
\end{tabular}
\caption{The number of tetrahedral elements and vertices for a sequence of hypersurface meshes for the rotating ellipsoid test case.} \label{ellipsoid_table}
\end{center}
\end{table}

\begin{figure}[h!]
\centering
\includegraphics[width=8cm]{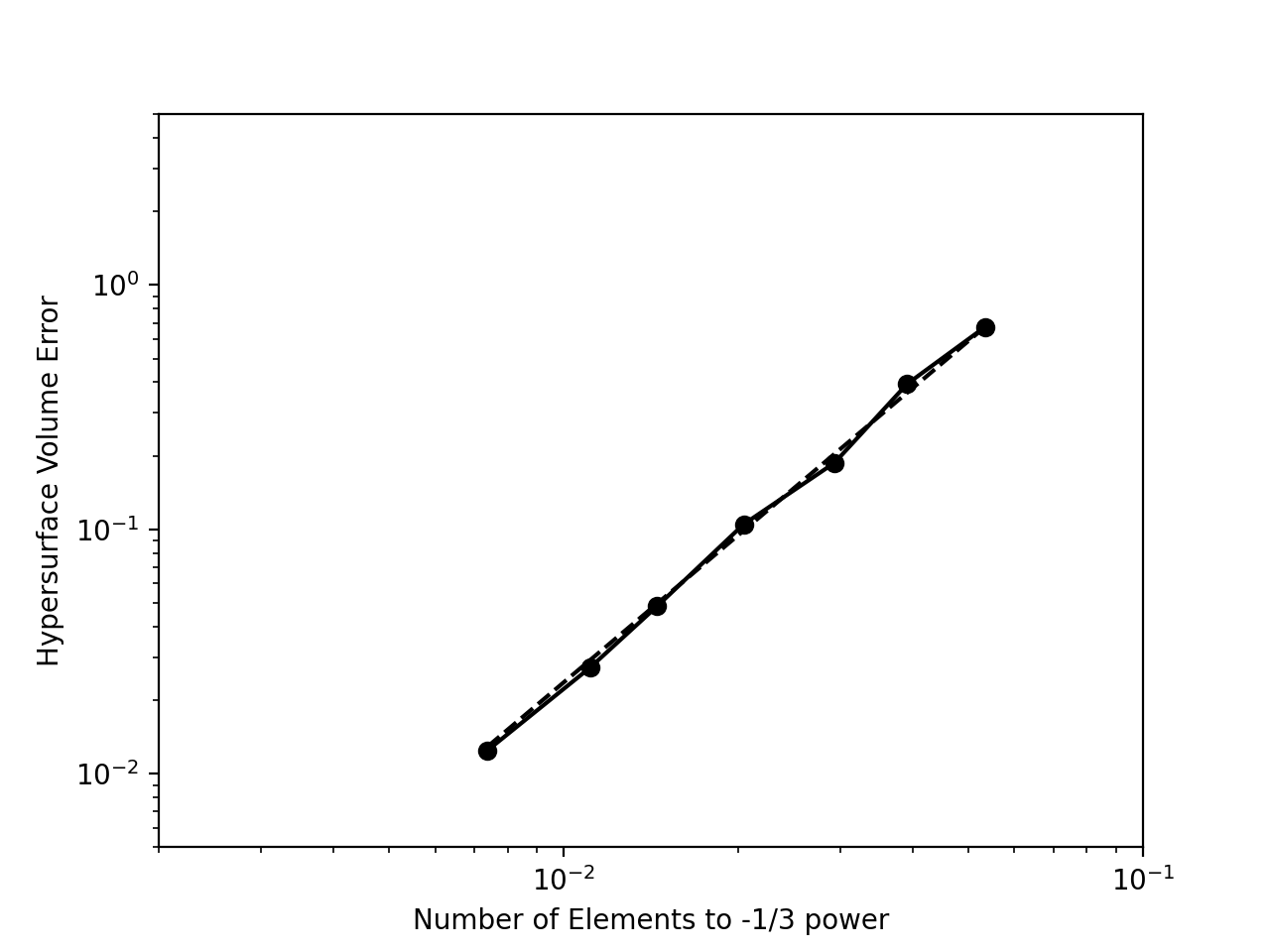}
\caption{Each point on the plot above represents the error between the volume of a hypersurface mesh and the exact volume for the rotating ellipsoid test case. The errors are plotted against the characteristic mesh spacing for a sequence of increasingly refined hypersurface meshes. In addition, a dashed line associated with 2nd-order convergence is plotted for reference.}
\label{rotating_ellipsoid_error_fig}
\end{figure}

We compared the volume of the hypersurface mesh at the final time ($t_f = 1$) against the exact volume of the hypersurface. The exact volume at $t_f = 1$ was calculated as follows
\begin{align*}
    V_{\text{exact}} (t_f) = L^3 - \frac{4}{3} \pi a b c.
\end{align*}
Figure~\ref{rotating_ellipsoid_error_fig} shows a plot of the volumetric error versus the mesh resolution. Second-order convergence is obtained, as expected.


\subsection{Rotating Tandem Ellipsoids}

\begin{figure}[h!]
\centering
\includegraphics[width=6cm,trim={2cm 0 2cm 0},clip]{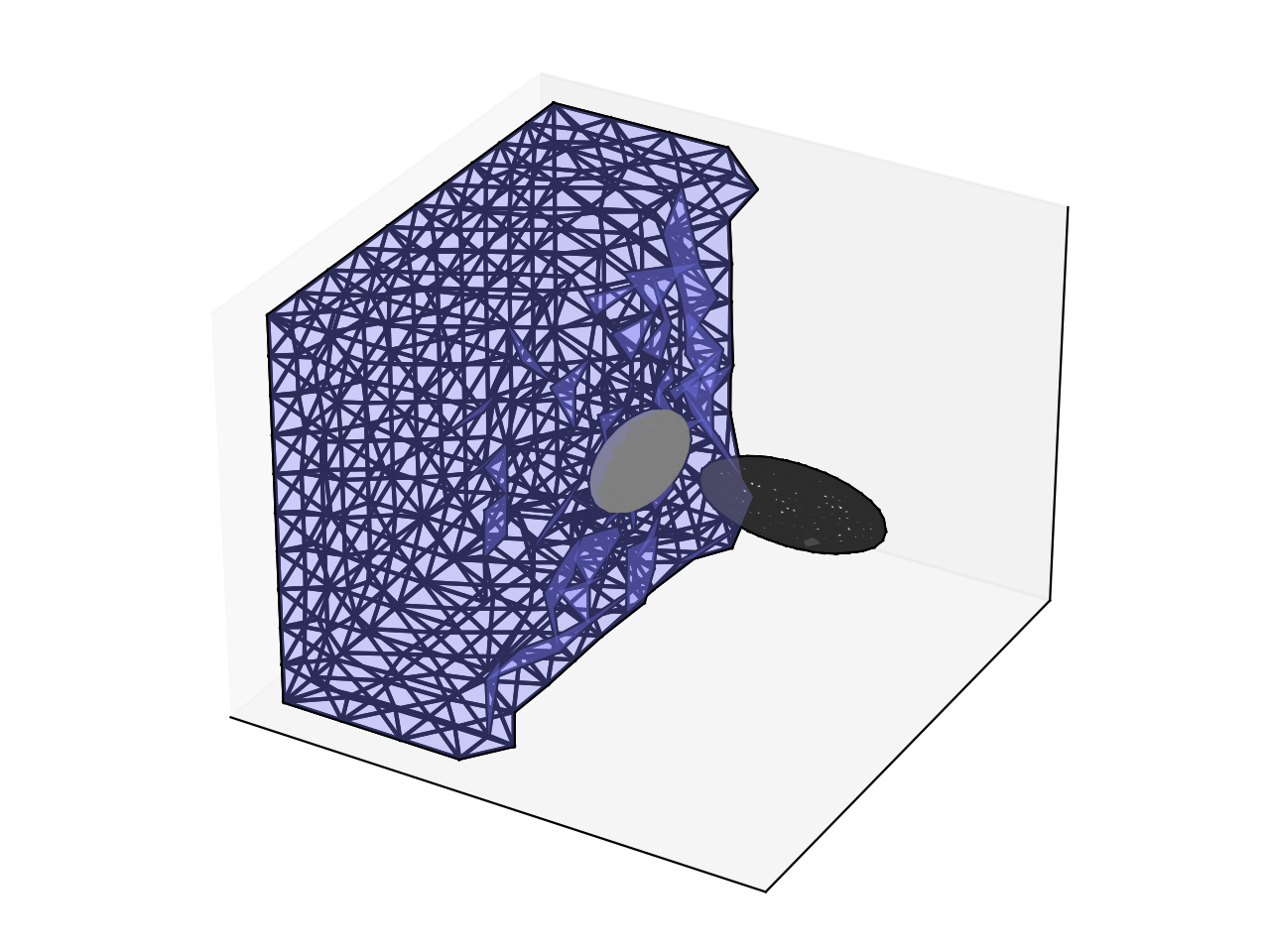} 
\includegraphics[width=6cm,trim={2cm 0 2cm 0},clip]{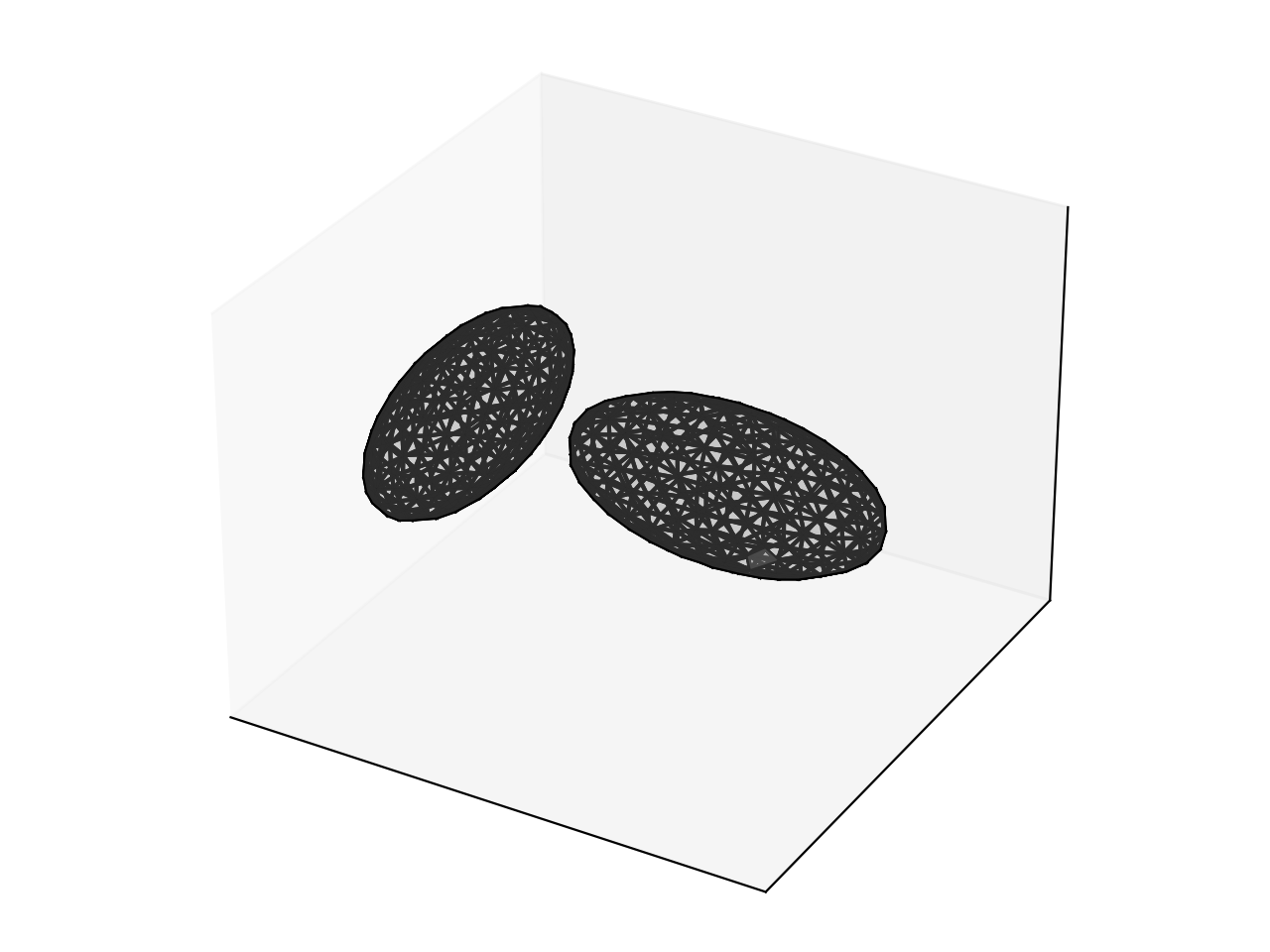}
\includegraphics[width=6cm,trim={2cm 0 2cm 0},clip]{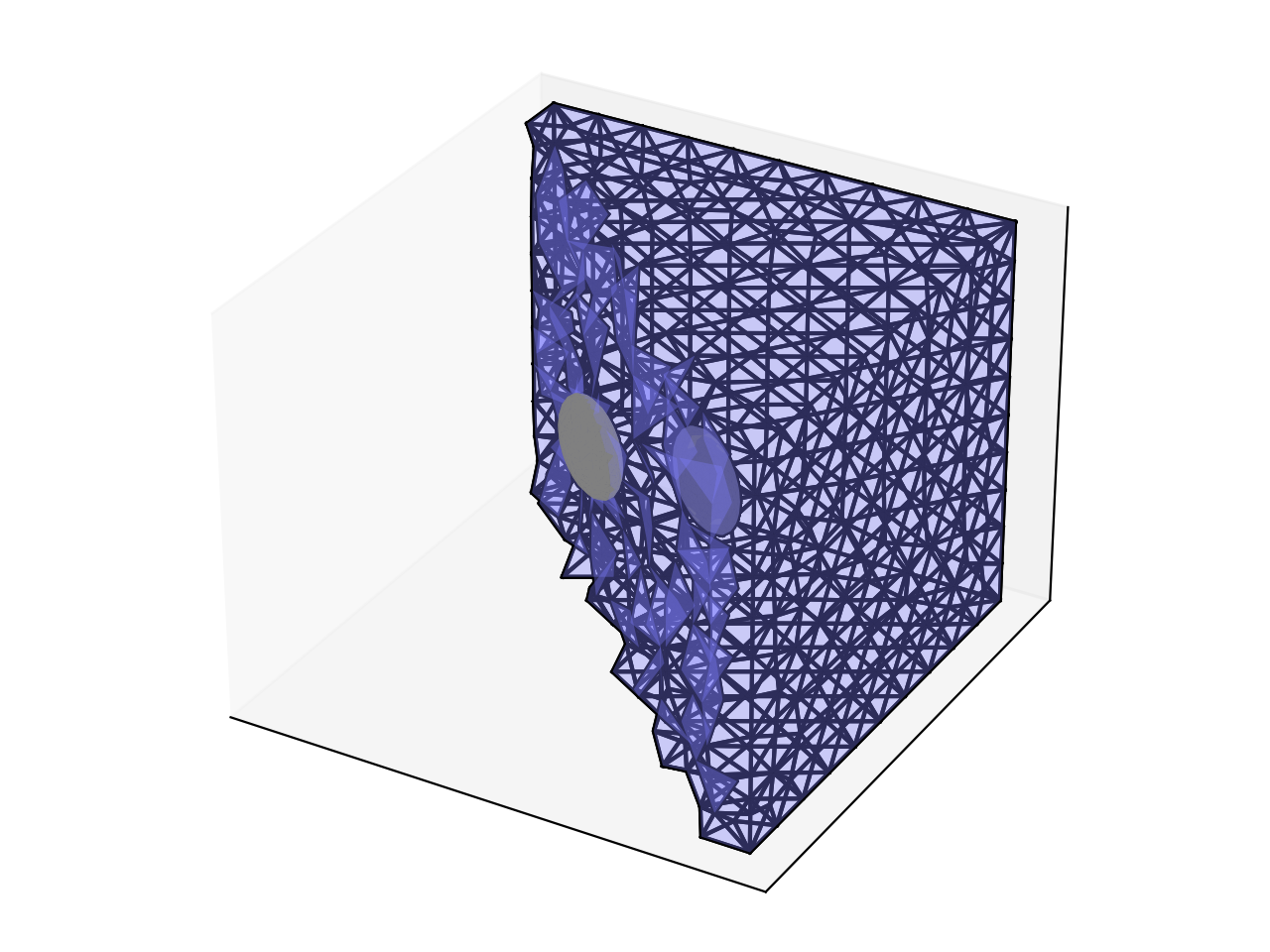} 
\includegraphics[width=6cm,trim={2cm 0 2cm 0},clip]{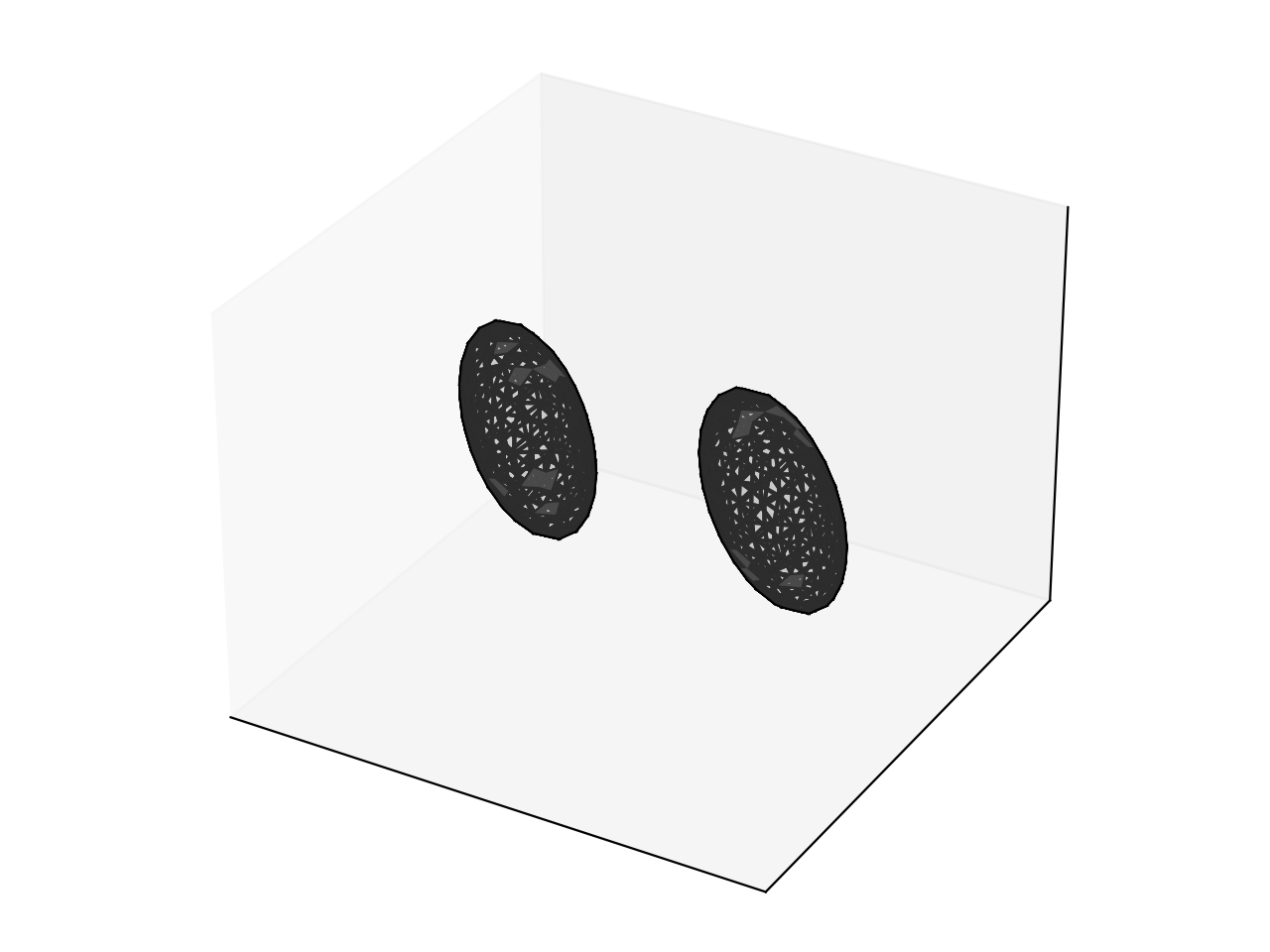}
\includegraphics[width=6cm,trim={2cm 0 2cm 0},clip]{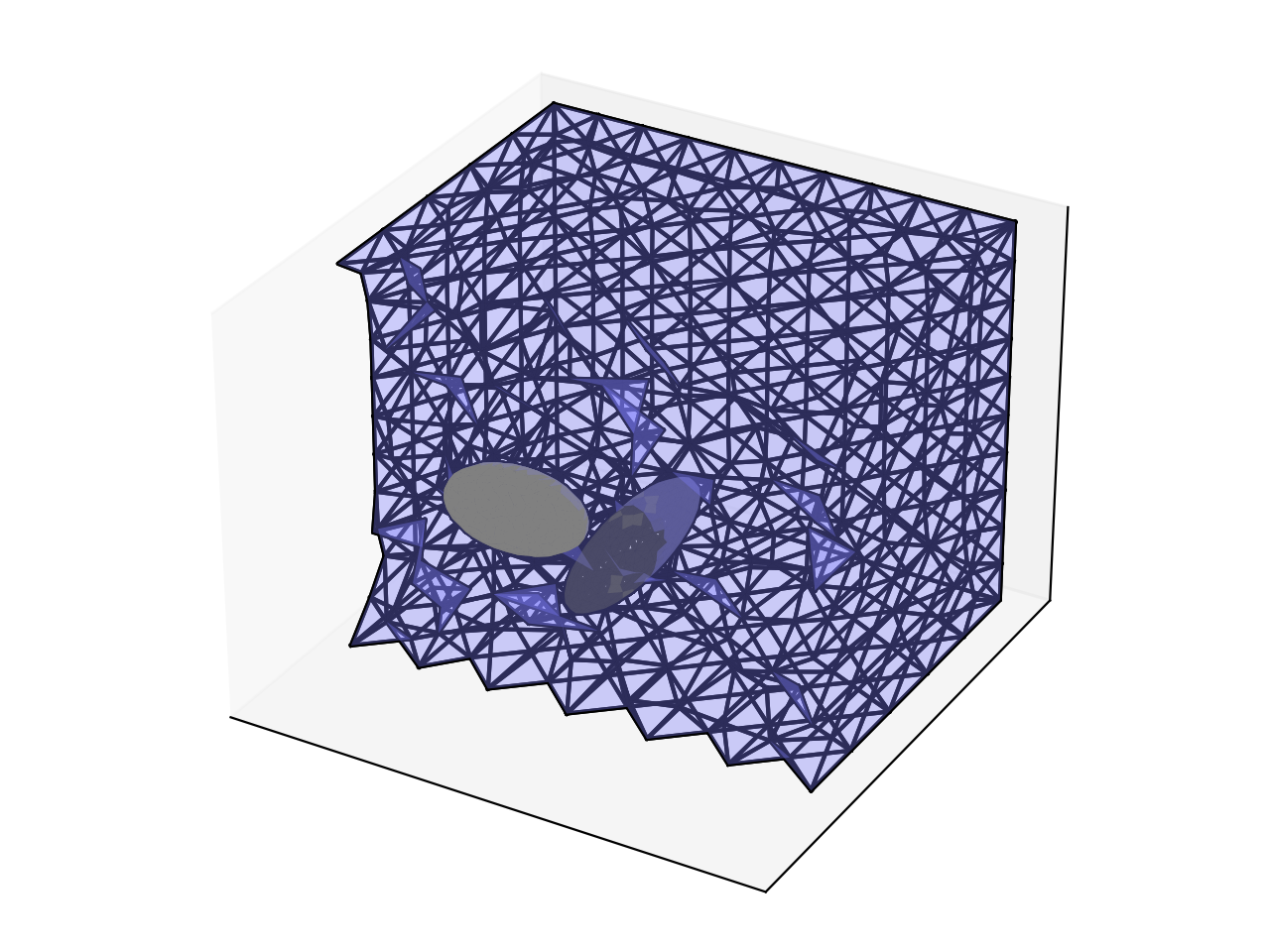} 
\includegraphics[width=6cm,trim={2cm 0 2cm 0},clip]{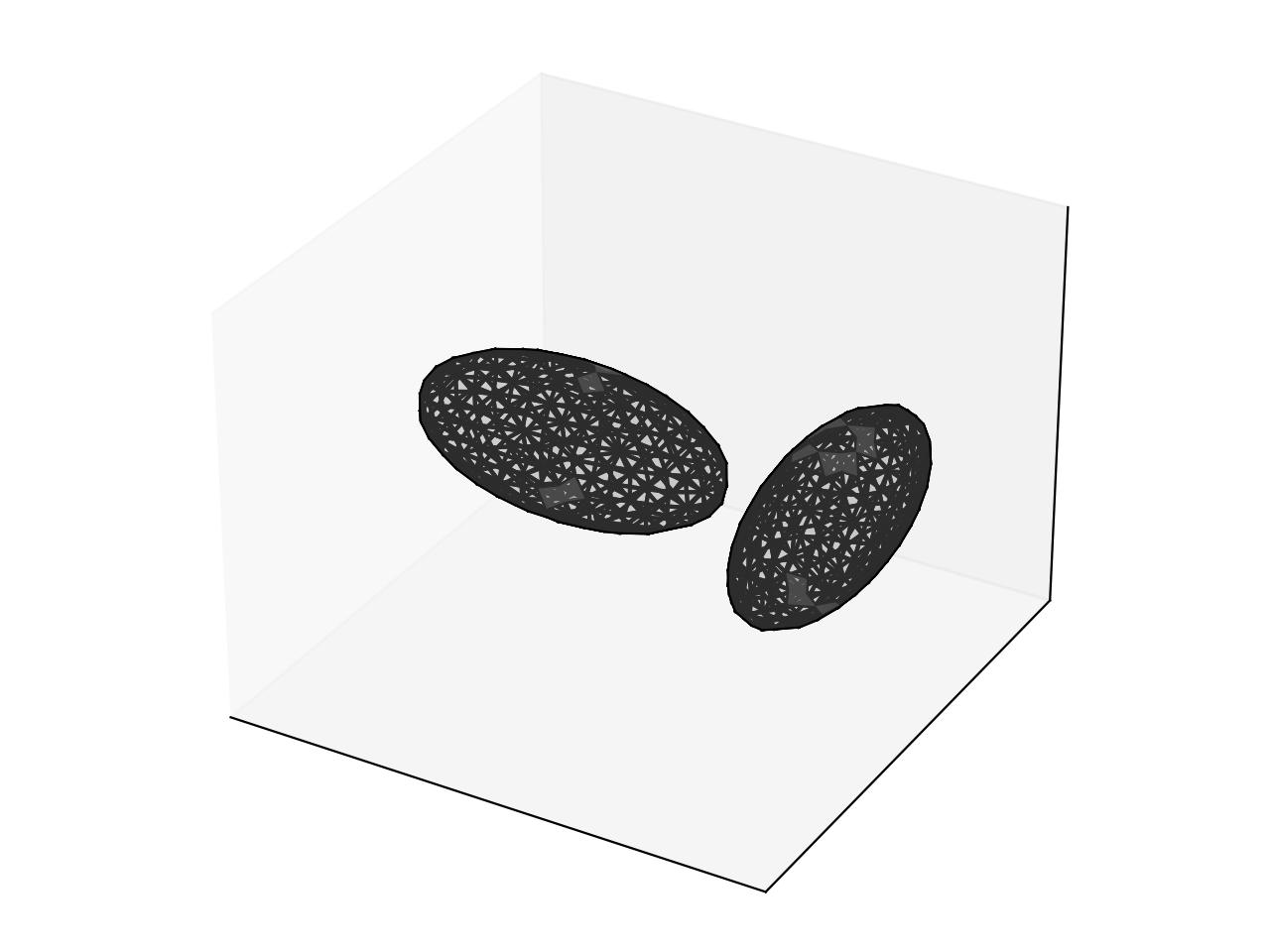}
\caption{Snapshots of rotating tandem ellipsoids at $t = 0$, $t = 0.5$, and $t = 1$, (top to bottom). On the left, a cross section of the coarsest hypersurface mesh is shown, alongside the ellipsoids' CAD. On the right, a zoomed-in view of the surface meshes on the ellipsoids is shown.}
\label{tandem_ellipsoids_snapshots}
\end{figure}

For this test case, the geometry consisted of a pair of ellipsoids with semi-axes of $a_1 = 1$, $b_1 = 3$, $c_1 = 2$ and $a_2 = 3$, $b_2 = 1$, $c_2 = 2$, respectively. Both ellipsoids were placed inside of a 3-cube with edge length $L = 20$, and bounds given by $\left[-7.5,12.5\right] \times \left[-10,10\right] \times \left[-10,10\right]$. The first ellipsoid was centered at $(0,0,0)$ and the second at $(5,0,0)$. In addition, the first ellipsoid rotated with angular velocity $(0, 0, \pi/2)$ rads/s, and the second rotated with velocity $(0, 0, -\pi/2)$ rad/s. On the time interval $[0,1]$, the motion of the ellipsoids created a pair of elliptical hyper-helixes that were contained inside of a tesseract. A family of hypersurface meshes was generated for this test case using the procedure described in the previous section. Table~\ref{tandem_ellipsoid_table} summarizes the properties of these meshes. Furthermore, Figure~\ref{tandem_ellipsoids_snapshots} shows several characteristic snapshots of the coarsest hypersurface mesh.   
\begin{table}[h!]
\begin{center}
\begin{tabular}{| c |r|r| }
\hline
Mesh & Elements & Vertices \\
\hline
1 & 603,432 & 128,055 \\
2 & 1,627,918 & 345,616 \\
3 & 4,523,078	& 959,759 \\
4 & 12,082,322 & 2,566,188 \\
5 & 35,245,617 & 7,478,525 \\
6 & 94,464,074 & 20,066,876 \\
7 & 278,561,736 & 59,127,488 \\
\hline
\end{tabular}
\caption{The number of tetrahedral elements and vertices for a sequence of hypersurface meshes for the rotating tandem ellipsoids test case.} \label{tandem_ellipsoid_table}
\end{center}
\end{table}

\begin{figure}[h!]
\centering
\includegraphics[width=8cm]{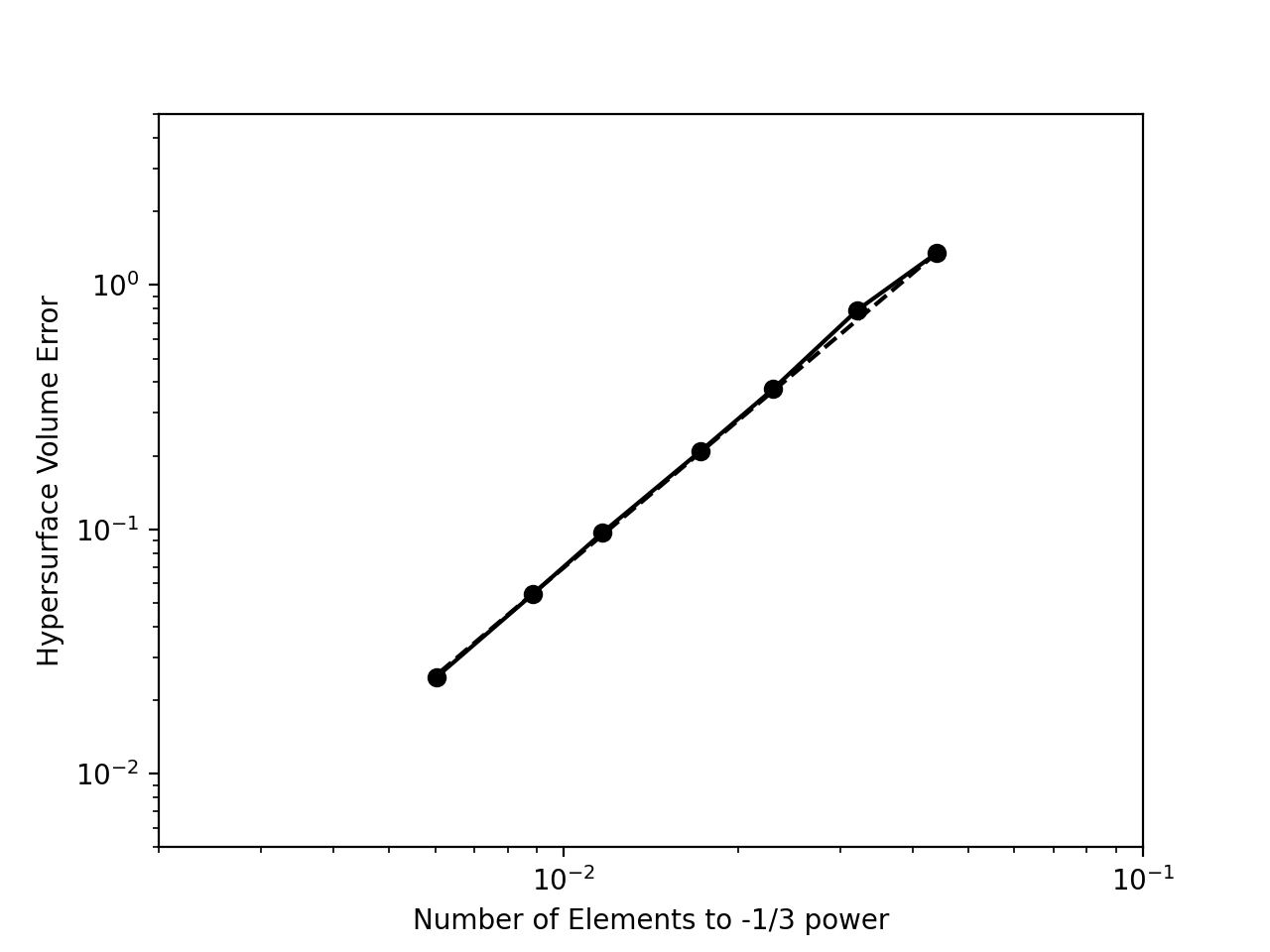}
\caption{Each point on the plot above represents the error between the volume of a hypersurface mesh and the exact volume for the rotating tandem ellipsoid test case. The errors are plotted against the characteristic mesh spacing for a sequence of increasingly refined hypersurface meshes. In addition, a dashed line associated with 2nd-order convergence is plotted for reference.}
\label{rotating_tandem_ellipsoid_error_fig}
\end{figure}

The exact hypersurface volume at final time $t_f = 1$ is given by
\begin{align*}
    V_{\text{exact}} (t_f) = L^3 - \frac{4}{3} \pi \left(a_1 b_1 c_1 + a_2 b_2 c_2\right).
\end{align*}
Figure~\ref{rotating_tandem_ellipsoid_error_fig} shows a plot of the volumetric error versus the mesh resolution. We obtain second-order convergence as expected. 


\section{Conclusion}
We have described in detail a general method for developing surface meshes in 2D+$t$ space time and hypersurface meshes in 3D+$t$ space time based on temporal planes (hyperplanes) derived from vertex trajectory-tracking through space time. These methods have been verified through numerical experiments by extruding/extending 2D and 3D objects along the temporal direction and comparing the approximate simplical surface areas or hypersurface volumes to the expected analytical results.  All numerical errors demonstrate 2nd-order convergence as the element densities of the surface (hypersurface) meshes increase, which demonstrates that our methods are working as expected.

In our future work, we will explore methods for Delaunay-based hypervolume meshing in 3D+$t$ space time. This work will include the development of methods for recovering a hypersurface boundary mesh once an initial 
 (unconstrained) hypervolume mesh has been generated.

\section*{Declaration of Competing Interests}

The authors declare that they have no known competing financial interests or personal relationships that could have appeared to influence the work reported in this paper.

\section*{Funding}

This research is sponsored by the Office of Naval Research, Code 351, through the Jet Noise Reduction Program under Program Officer, Dr. Steven Martens. Penn State University received funding under contract number N00173-22-2-C008.

\pagebreak
\clearpage

{\footnotesize\bibliography{technical-refs}}

\end{document}